\documentclass[journal]{aiaa-pretty}    % journal, submit, article
\usepackage[utf8]{inputenc}  % Helps a lot with arXiv submission!!!
\usepackage{amssymb}
\usepackage{amsmath}
\usepackage{booktabs}
 \usepackage{varioref}%  smart page, figure, table, and equation referencing
 \usepackage{wrapfig}%   wrap figures/tables in text (i.e., Di Vinci style)
 \usepackage{threeparttable}% tables with footnotes
 \usepackage{dcolumn}%   decimal-aligned tabular math columns
  \newcolumntype{d}{D{.}{.}{-1}}
 \usepackage{nomencl}%   nomenclature generation via makeindex
  \makeglossary
 \usepackage{subfigure}% subcaptions for subfigures
 \usepackage{graphicx}
 \usepackage{subfigmat}% matrices of similar subfigures, aka small mulitples
 \usepackage{fancyvrb}%  extended verbatim environments
  \fvset{fontsize=\footnotesize,xleftmargin=2em}
 \usepackage{lettrine}%  dropped capital letter at beginning of paragraph
\usepackage{setspace}
\usepackage{cite}

\newtheorem{theorem}{Theorem}
\newtheorem{proposition}{Proposition}
\newtheorem{remark}{Remark}
\newtheorem{lemma}{Lemma}
\newtheorem{definition}{Definition}

\newtheorem{hypothesis}{Hypothesis}
\newtheorem{corollary}{Corollary}
\newenvironment{proof}{{\it Proof. }}{\hfill $\Box$}

\newcommand{\Real}{\mathbb R}
\newcommand{\set}[1]{\left\{#1\right\}}
\newcommand{\real}[1]{{\mathbb R}^{#1}}
\newcommand{\bb}{{\boldsymbol b}}

\newcommand{\bff}{{\boldsymbol f}}

\newcommand{\bu}{{\boldsymbol u}}

\newcommand{\bw}{{\boldsymbol w}}
\newcommand{\bx}{\boldsymbol x}
\newcommand{\by}{{\boldsymbol y}}

\newcommand{\bB}{{\boldsymbol B}}

\newcommand{\bK}{{\boldsymbol K}}

\newcommand{\bdelta}{\mbox{\boldmath$\delta$}}

\author{ %
I. M. Ross\thanks{Distinguished Professor and Program Director, Control and Optimization, Department of Mechanical and Aerospace Engineering.\texttt{imross@nps.edu}}
\\
\textit{Naval Postgraduate School, Monterey, CA 93943}
}
\title{A Universal Birkhoff Theory for Fast Trajectory Optimization}

\abstract{
Over the last two decades,  pseudospectral methods based on Lagrange interpolants have flourished in solving trajectory optimization problems and their flight implementations.  In a seemingly unjustified departure from these successes, a new starting point for trajectory optimization is proposed.  This starting point is based on the recently-developed concept of universal Birkhoff interpolants.  The  new first-principles' approach offers a substantial computational upgrade to the Lagrange theory in completely flattening the rapid growth of the condition numbers from $\mathcal{O}(N^2)$ to $\mathcal{O}(1)$, where $N$ is the number of grid points. Furthermore, a Birkhoff discretization offers the theoretical possibility of an infinite-order rate of convergence. In addition, the Birkhoff-specific primal-dual computations are isolated to a well-conditioned linear/convex system even for nonlinear, nonconvex problems. This is Part~I of a two-part paper.  In Part~I, the new theory is developed on the basis of two hypotheses. Other than these hypotheses, the theory makes no assumptions on the choices of basis functions, be it polynomials or otherwise, or even the distribution of grid points.
Several theorems are proved to establish the mathematical equivalence between direct and indirect Birkhoff methods.  In Part~II of this paper, it is shown that a select family of Gegenbauer grids satisfy the two hypotheses required for the theory to hold.  The theory developed in Part~I is used in Part~II to produce the fast computational implementation.
}
\begin{document}
\maketitle

%\end{spacing}

%\newpage
\section{Introduction}\label{sec:intro}
\begin{figure}[h!]
      \centering
      {\parbox{\columnwidth}{
      \centering
      {\includegraphics[width = 0.95\columnwidth]{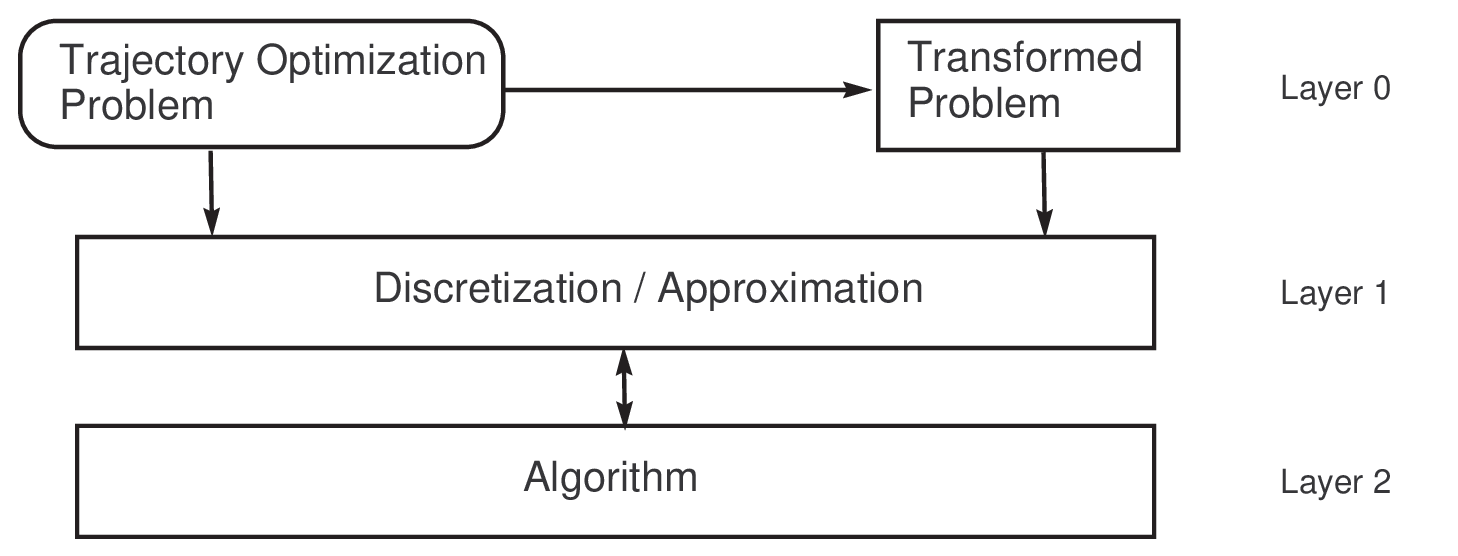}}
      \caption{{A schematic of three layers that captures the essence of diverse trajectory optimization techniques.}}
      \label{fig:TrajOptPerspective}
      }
      }
\end{figure}

Barring a few exceptions, trajectory optimization problems are open-loop optimal control problems\cite{longuski,ross-book}. Purportedly, there is a wide variety of methods to solve such problems\cite{conway:survey,longuski,ross-book,perspective,trelat-2017, trelat:survey,trelat:guidance, lu:entry-guidance, asteroid-trajOpt,ddp-2021,seywald-94,vonStryk-92,convex-2015,convex-2011,flatness-2013} that range from application-specific techniques\cite{longuski,lu:entry-guidance,trelat-2017,trelat:guidance,asteroid-trajOpt,trelat:survey,convex-2015} to general-purpose methods\cite{ross-book,longuski,perspective,conway:survey,trelat:survey,vonStryk-92}.  Moreover, within each of these categories, there are further ideas that range from semi-analytic or geometric methods\cite{trelat:guidance,perspective,longuski,lu:entry-guidance,convex-2015} to various subclasses of generic problems such as linear, convex, differentially-flat, unconstrained etc.\cite{trelat:survey,ddp-2021,convex-2015,convex-2011, seywald-94,perspective,flatness-2013}.  One possible way to view these diverse set of ideas is through a prism of three layers identified in \cite{perspective}: (1) Transformations, (2) Discretizations (or approximations) and (3) Algorithms. From this perspective, a classical indirect method\cite{vonStryk-92,trelat:survey} may be viewed as first transforming a given problem to its primal-dual form (see Layer 0 in Fig.~\ref{fig:TrajOptPerspective}) and subsequently applying a discretization (Layer 1) and an algorithm (Layer 2) to solve the transformed, generalized boundary-value problem (G-BVP)\cite{ross-book,longuski,trelat-2017}.
%======================================================================================
%==========================================================================================
%
In a flatness-based approach\cite{flatness-2013,perspective}, the trajectory optimization problem is first transformed to a differentially-flat output space (as part of Layer~0 of Fig.~\ref{fig:TrajOptPerspective}) before a discretization or approximation is performed. The flatness-based approach is a direct method; however, a classical direct method\cite{conway:survey,vonStryk-92} largely involves Layers~1 and 2, where, the second layer is frequently implemented using established optimization routines. Special-purpose techniques\cite{convex-2011,lu:entry-guidance,asteroid-trajOpt,trelat:guidance} involve designing or using specific algorithms (Layer~2 in Fig~\ref{fig:TrajOptPerspective}) that are suited for a particular problem or a specific-class of problems (see Layer~0 in Fig.~\ref{fig:TrajOptPerspective}).
Depending upon the transformation employed, Layer~0 can be labor-intensive\cite{longuski,vonStryk-92} or problem-specific\cite{seywald-94,convex-2015} or automated\cite{ross-book,trelat:survey}. See \cite{perspective} for further details on this perspective.  In any event, no matter what ideas are used, the ultimate goal of all methods is to solve a given problem in the most efficient manner.  One of the metrics of efficiency is computational speed. Deferring a discussion of the meaning of computational speed, a key driver for fast trajectory optimization, particularly in aerospace application, is closed-loop guidance\cite{ross-book,lu:editorial}. If a solution to a trajectory optimization problem can be computed in ``real time,'' then optimal feedback control is possible\cite{ross-book,RTOC-jgcd,lu:entry-guidance,mease:guidance-compute,trelat:guidance,HJB=PP=DIDO} without a need to solve the difficult Hamilton-Jacobi-Bellman (HJB) equation\cite{brysonHo,vinter,clarke-2013book}.  But what exactly is real time?  Intuitively, it is apparent that if a dynamical system is ``slow'' then the time to recompute an optimal open-loop control can be large. Conversely, a faster dynamical system necessitates a faster computational speed to achieve the same closed-loop-objectives of a slower system.  This intuitive concept is made rigorous in \cite{RTOC-jgcd} and \cite{LGR-rtoc} in terms of $\text{Lip}_x\bff $, the Lipschitz constant\cite{clarke-2013book,vinter} with respect to $\bx$ of the dynamical system, $\dot\bx = \bff(\bx, \bu, t) $ (where, the symbols have the usual signficance\cite{ross-book,brysonHo,vinter}).  That is, if a solution to a trajectory optimization problem can be recomputed within a time interval $\tau_c$ given by the inverse of the Lipschitz constant,
\begin{equation}\label{eq:RTOC-def}
\tau_c = \frac{1}{\text{Lip}_x\bff }
\end{equation}
then the resulting algorithm is termed real-time\cite{ross-book,RTOC-jgcd,LGR-rtoc}.
This foundational concept of real-time is agnostic to computational hardware; in fact, it levies the requirements for pairing the specifications of a computing system to a specific dynamical system.  On the basis of \eqref{eq:RTOC-def}, real-time optimal control (RTOC) was used in \cite{LGR-rtoc} to stabilize an otherwise-unstable inverted pendulum without explicitly constructing a (nonlinear) feedback controller, $\bu = \bK(\bx,t)$ \cite{brysonHo,vinter,clarke-2013book}.  Equation~\eqref{eq:RTOC-def} was also used in \cite{RTOC-jgcd} to ground-test various optimal feedback slew maneuvers for NPSAT1, an experimental satellite built at the Naval Postgraduate School and launched in June~2019.  These engineering applications demonstrated the practical use of the theoretical concept of using the Lipschitz frequency to define RTOC.  Equation~\eqref{eq:RTOC-def} facilitates the construction of a \emph{closed-loop} control system without the requirement for producing a \emph{closed-form} solution to a feedback control law\cite{ross-book}.  Because a traditional closed-form feedback controller automatically implies real-time computation of the function, $\bK(\cdot)$, it is apparent that real-time computation is more fundamental to the production of closed-loop controllers than ``analytical expressions'' of $\bK(\bx,t)$ in terms of state variables.  Note that RTOC does not mean traditional closed-form feedback controls must be abandoned. In fact, as described in \cite{ross-book}, RTOC can be used in concert with heritage control systems as an enabling technology kernel for designing various outer-loops under the rubric of different guidance and control system architectures\cite{lu:editorial}.

It is important to recognize that application-specific RTOC algorithms are not new; they have been used for quite some time in the aerospace industry under the heading of closed-loop guidance\cite{ross-book,lu:prop-guidance,PEG:history,mease:guidance-compute,lu:entry-guidance}.  For example, the well-known powered explicit guidance\cite{PEG:history} is a closed-loop application of the optimal open-loop linear tangent equation\cite{brysonHo,ross-book}. Many more application-specific RTOC algorithms have been developed in recent years\cite{lu:prop-guidance,mease:guidance-compute,lu:entry-guidance}. General-purpose RTOC-based feedback controllers are, relatively, more recent developments\cite{ross-book, RTOC-jgcd,HJB=PP=DIDO}. As shown in\cite{RTOC-jgcd,HJB=PP=DIDO, bollinoRossDoman-2006,bollinoRoss-2007, lewis-NATO-2007, hurniSekRoss-2009, bollinoLewis-2008, bollinoLewis-2007,bollinoEtAl-2007}, the same general-purpose trajectory optimization method used to stabilize an inverted pendulum in \cite{LGR-rtoc} was reused multiple times, without change, to generate attitude guidance commands for NPSAT1\cite{RTOC-jgcd}, entry guidance algorithms for hypersonic vehicles\cite{bollinoRossDoman-2006,bollinoRoss-2007}, autonomous operations of ground robots in cluttered and uncertain environments\cite{lewis-NATO-2007,hurniSekRoss-2009,HJB=PP=DIDO}, cooperative controls for multiple uninhabited aerial vehicles\cite{bollinoLewis-2008}, and other disparate applications\cite{bollinoLewis-2007,bollinoEtAl-2007}.  The common theme in all these diverse applications was pseudospectral (PS) optimal control theory\cite{PSReview-ARC-2012}.
Unlike a classical direct or indirect method\cite{vonStryk-92}, PS optimal control theory is based on the inherent coupling of all elements of the trajectory optimization layers indicated in Fig.~\ref{fig:TrajOptPerspective}\cite{spec-alg,PSReview-ARC-2012, DIDO:arXiv}.  Note that simply using a PS discretization as a classical direct method does not qualify the approach to be full-fledged member of PS optimal control theory.   Briefly put, a primal PS discretization is first chosen\cite{PSReview-ARC-2012} by levying the requirement that it commute with its dual transformation\cite{ross:roadmap-2005} (Layers~0 and 1).  This commutative property is thereafter used to generate a guess-free spectral algorithm\cite{spec-alg, DIDO:arXiv,ross:guess-free} (Layer~2) that interacts with the primal/dual PS discretization  and vice versa (Layer~1) to solve the transformed problem (Layer~0) in a manner that is consistent with an application of Pontryagin's Principle\cite{ross-book} to the (original) trajectory optimization problem indicated in Layer~0.  See \cite{DIDO:arXiv} and the references contained therein for details.  For all these reasons, a PS-control-theoretic implementation cannot be classified as either a direct or indirect method but is explainable in terms of the fundamental constituent elements of trajectory optimization indicated in Fig.~\ref{fig:TrajOptPerspective}.

\begin{figure}[h!]
      \centering
      {\parbox{\columnwidth}{
      \centering
      {\includegraphics[width = 0.95\columnwidth]{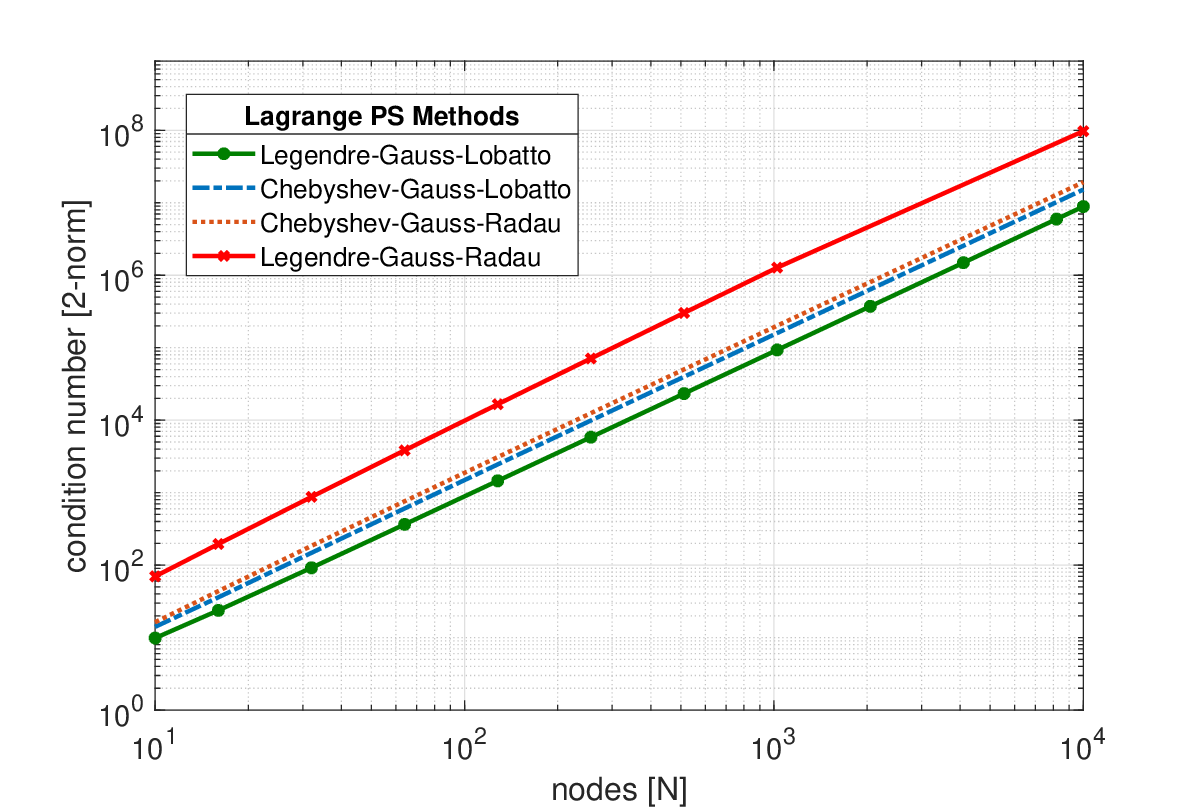}}
      \caption{{$\mathcal{O}(N^2)$ growth in the condition number of a Lagrange PS method for various choices of grid selections; data from \cite{fastmesh}.}}\label{fig:condNums4Lags}
      }
      }
\end{figure}

One of the primary reasons for using PS optimal control theory for RTOC (as a general-purpose approach) is its near-exponential rate of convergence\cite{Kang_2008_convergence,kang-rate-2010} (with respect to Layers~0 and 1 in Fig.~\ref{fig:TrajOptPerspective}).  That is, it takes relatively few grid points to generate an accurate trajectory solution and even fewer points to close the loop\cite{RTOC-jgcd,LGR-rtoc,bellman-conf,bellman-low-t}.  The latter result is due to the anti-aliasing property of the ensuing Bellman PS method\cite{bellman-low-t, bellman-conf,TEI-JGCD-2011}.  For problems where there is a need for a highly dense grid, the theoretical fast rate of convergence of a PS discretization is practically dampened by round-off errors\cite{twist,trefethen-2000} and the rapidly increasing condition number of the resulting linear algebra that underpins the trajectory optimization algorithm\cite{spec-alg,ross:guess-free,rossJCAM-1,DIDO:arXiv} (see Layer~2 in Fig.~\ref{fig:TrajOptPerspective}).  The latter point is illustrated in Fig.~\ref{fig:condNums4Lags} which shows the $\mathcal{O}(N^2)$ growth in the condition number of a PS discretization, where $N$ is the number of grid points.
%======================================================================================
%==========================================================================================
High condition numbers are anathema to both computational speed and accuracy\cite{trefethen-bau-1997}.
To mitigate the impact of high condition numbers, preconditioners are commonly used\cite{hesthavan,elbarbary} to ``bend the curve'' of Fig.~\ref{fig:condNums4Lags}. Prior to the work of Wang et al\cite{wang}, the best preconditioner generated $\mathcal{O}\left(\sqrt{N}\right)$ condition numbers (caveated by certain technical conditions\cite{hesthavan,elbarbary,wang}). As impressive as $\mathcal{O}\left(\sqrt{N}\right)$ is, the results of \cite{wang} changed everything.  This is because:
\begin{enumerate}
\item It solved the long-standing problem\cite{hesthavan,elbarbary,wang} of finding the computationally perfect preconditioners in generating $\mathcal{O}(1)$ condition numbers; and,
\item It opened the door for abandoning the use of Lagrange interpolants as a starting point for all PS methods\cite{fornberg,boyd,trefethen-2000}; i.e., it effectively made all prior PS methods ``obsolete.''
\end{enumerate}
As a consequence of the first aforementioned point, fast and accurate solutions to BVPs became possible using the perfect preconditioner.  Nonetheless, it is the second point that made the work of Wang et al\cite{wang} quite revolutionary.  This is because they showed how a PS method could be conceived from a new starting point. That is, instead of using Lagrange interpolants and subsequent preconditioners to formulate a fast PS method, it was shown in \cite{wang} that the better starting point was a Birkhoff interpolant\cite{lorentz,schoenberg}. This essentially rendered all prior methods obsolete.  Thus, if a PS method were to be introduced using Birkhoff interpolants, then many of the challenges associated with Lagrange interpolants vanish in one fell swoop. Furthermore, because low condition numbers enhance both computational speed and accuracy\cite{trefethen-bau-1997}, the resulting well-conditioned Birkhoff PS methods are capable of generating fast and accurate solutions over thousands of grid points\cite{wang,fastmesh}.  Thus, a Birkhoff PS method, if applied to trajectory optimization, offers itself as a prime practical candidate for a general-purpose RTOC method to solve nonlinear, nonconvex aerospace guidance problems.  This paper advances the foundations of such methods through the use of the new universal Birkhoff interpolants proposed in \cite{newBirk-2023}. % Ross et al.
These universal Birkhoff interpolants are different from those of Wang el al\cite{wang} in several different ways\cite{newBirk-2023} and offer many theoretical and practical advantages, particularly in dual space transformations (see Layer~0 in Fig.~\ref{fig:TrajOptPerspective}).

A Birkhoff interpolant\cite{lorentz,schoenberg} allows a heterogenous mix of derivatives of various orders as part of its constituent elements. A careful selection of this mixture can be used to generate a family of fast and accurate methods for solving a variety of BVPs\cite{wang}.  Motivated by these ideas, we adapted the results of \cite{wang} for solving trajectory optimization problems\cite{fastmesh,furtherResults}. Although trajectory optimization problems may be viewed as generators of BVPs, an adaption of \cite{wang} to solve optimal control problems proved to be not entirely straightforward\cite{fastmesh,furtherResults}.  This is because \cite{wang} requires a production of different Birkhoff interpolants for different BVPs based on their order and the type of boundary conditions such as Dirichlet, Neumann or Robin. In seeking a more universal approach, we advanced in \cite{newBirk-2023}
a new family of Birkhoff interpolants that eliminate the need for such specialized treatments.
These universal Birkhoff interpolants of \cite{newBirk-2023} address BVPs of various orders and mixtures of boundary conditions without the need to redevelop the method for different problems. In this two-part paper, we further the results of \cite{newBirk-2023} to develop a new approach for solving general-purpose trajectory optimization problems.  Part II of this paper is \cite{newBirk-part-II}.

In this paper (i.e., Part I), we show that the new Birkhoff interpolants satisfy the covector mapping principle (CMP)\cite{ross-book, perspective, ross:roadmap-2005,CMP-CDC-2006} under certain hypotheses pertaining to the selection of grid points.  These hypotheses are shown to be satisfied in Part II of this paper when the grid is selected from a family of Gegenbauer node points.  A satisfaction of the CMP by any discretization method has major ramifications on trajectory optimization.  This is because the CMP implies the following\cite{ross-book, perspective,spec-alg, DIDO:arXiv,scaling}:
\begin{enumerate}
\item The oft quoted\cite{conway:survey,vonStryk-92,cots-2017} distinctions between a direct and indirect method vanishes\cite{ross-book,perspective};
\item An equivalence between a direct and indirect method can be advantageously used for the production of significantly faster algorithms\cite{spec-alg, DIDO:arXiv,rossJCAM-1,ross-CD};
\item Every computed solution for a specific problem can be readily tested for optimality by simple checks of key necessary conditions resulting from an application of Pontryagin's Principle\cite{ross-book,DIDO:arXiv,newBirk-part-II};
\item If a problem is badly scaled because of the use of engineering units, it can be quickly re-scaled for computational speed and accuracy by employing the concept of (non-canonical) designer units\cite{ross-book, scaling} made possible by covector balancing techniques\cite{scaling}; and,
\item The connection between the HJB equations and Pontryagin's Principle\cite{vinter,brysonHo,clarke-2013book} can be rigorously utilized for optimality proofs and practical validation of the computed solution even for nonconvex problems\cite{HJB=PP=DIDO}.
\end{enumerate}
Taken all together, the promise of a universal Birkhoff theory for trajectory optimization has significant implications across all of its three layers depicted in Fig.~\ref{fig:TrajOptPerspective}. Its applications to RTOC are as follows:
%
%==================================
\begin{enumerate}
\item Theoretically, a Birkhoff interpolant offers the possibility of an infinite-order rate of convergence (see Section~\ref{sec:fast=compcomp}). This capability implies the resulting scale of the Birkhoff-discretized mathematical programming problem can be ``exponentially smaller'' when compared to other methods\cite{TAC:linearizable,paris-comparision-2006};
\item The growth in the condition numbers of a key matrix block of the linear algebra that underpins the accuracy of the resulting trajectory optimization algorithm is flattened to $\mathcal{O}(1)$;
%
%======================================================================================
\begin{figure}[h!]
      \centering
      {\parbox{0.9\columnwidth}{
      \centering
      {\includegraphics[width = 0.95\columnwidth]{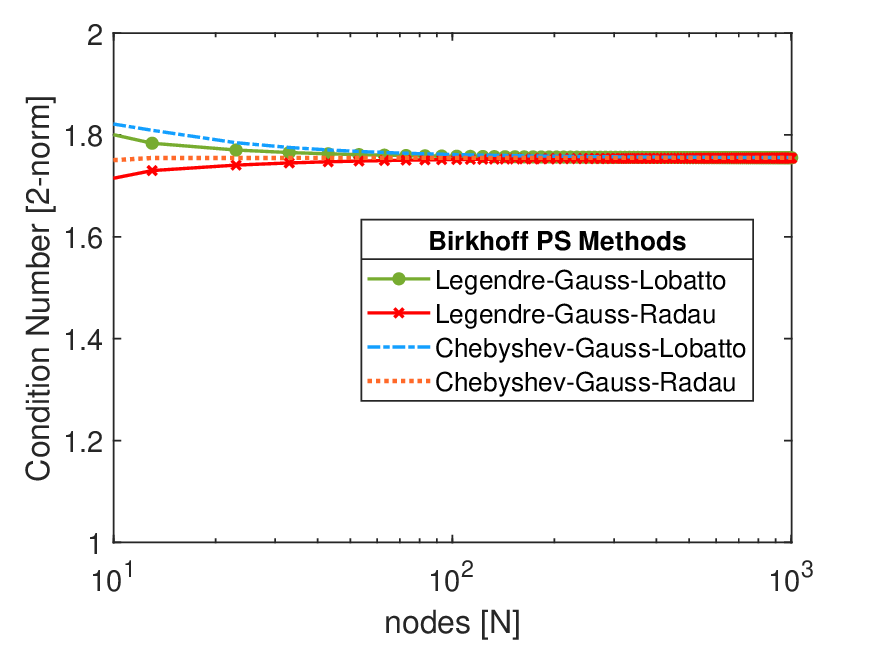}}
      \caption{\textsf{$\mathcal{O}(1)$ variation in the condition number of a discretization based on using the universal Birkhoff interpolants\cite{newBirk-2023} for various choices of grid selections.}}\label{fig:condNums4Birks}
      }
      }
\end{figure}
%==========================================================================================
%
see Fig.~\ref{fig:condNums4Birks} and  \cite{newBirk-part-II}.
This fact implies one can generate fast and accurate solutions;
\item Round-off errors\cite{twist,wang} that plague the production of large differentiation matrices vanish in one fell swoop.  This is because the Birkhoff matrices can be computed over thousands of grid points in a fast and stable manner\cite{wang,newBirk-2023,sandia-aas-23};
\item Coupled with recent advancements in spectral methods\cite{Bogaert-2014,Hale-fast-comp-2013,olver-2021,Gil-2019}, even millions of grid points can be generated at an $\mathcal{O}(1)$ computational speed. This fact implies grid refinement can be performed ``instantaneously\cite{fastmesh}.''
\item If one chooses a Chebyshev grid, the Birkhoff matrix-vector products can be computed at an $\mathcal{O}(N\log(N))$ computational speed\cite{sandia-aas-23} using a fast Fourier transform (FFT)\cite{kreyszig-2011,fft-original,fftw}.
\end{enumerate}
%==================================
%

\section{What is Fast Trajectory Optimization?}\label{sec:fast=compcomp}
Before delving into the details of the universal Birkhoff theory and its computational speed, it is important to delineate a more precise meaning of the word, fast.  From a simple-minded perspective, one can just define fast as computer time.  This definition is rife with naivety.  One obvious drawback of this definition is that a trajectory optimization software can be made to run faster by simply implementing it in the form of a compiled code (e.g., C/C++) versus an interpreted program (e.g., MATLAB${}^{TM}$).  Furthermore, within a given computer implementation, a trajectory optimization code can be made faster by providing Jacobian/Hessian information. A number of other factors, such as initial guess, accuracy, coordinate transformations (see Layer~0 in Fig.~\ref{fig:TrajOptPerspective}) etc. can affect run time. Even worse, human factors, such as bad coding, can affect run times in a deleterious manner.  In view of all such distortions, and much more, we do not use computer time as a metric for computational speed.

In this paper, we use the word fast in terms of the fundamental mathematics of computational complexity\cite{luenberger-2008,nesterov-book-2004,NW:NumOptBook,bazaraa-2006} that is invariant to all of the distortions noted in the previous paragraph.  Thus, for example, a computational process that is $\mathcal{O}(n\log(n))$ is faster than $\mathcal{O}(n^2)$ regardless of whether it is implemented in a complied form or an interpreted program or by an inexperienced coder.  In the previous sentence, $n$ is a generic integer that defines some scale. For the purposes of completeness, we note that such measures of computational speed are typically applicable asymptotically and may not hold for small values of $n$; see, for example, Fig.~\ref{fig:condNums4Birks}.  Furthermore, because the notation $\mathcal{O}(\cdot)$ hides a constant, it is quite possible that a $\mathcal{O}(n^2)$ operation may run faster than $\mathcal{O}(n\log(n))$ for sufficiently small values of $n$.  In any case, to the best of the author's knowledge, there is currently no comprehensive theory of computational complexity for trajectory optimization. Short of constructing one in this paper, we use the computational complexity of the components of trajectory optimization, delineated in Fig.~\ref{fig:TrajOptPerspective}, to define computational speed with appropriate caveats.  In this paper, we focus our attention on the components of Layer~1 and its strong dependence on Layer~2\cite{perspective,spec-alg,DIDO:arXiv,rossJCAM-1,ross-CD} (and vice versa) as indicated by the double-headed arrow in Fig.~\ref{fig:TrajOptPerspective}.

\subsection{Foundations of Infinite-Order Convergence Rate}

A central feature of computational speed pertaining to Layer~1 is the rate of convergence of the discretization\cite{perspective,DIDO:arXiv}.  A fast rate of convergence is desired for fast trajectory optimization.  The theoretical rates of convergence from Euler to PS discretizations are shown in Table~\ref{table:convergence-rates}.
%-------------------- TABLE -------------------------------
%
\begin{table*}[t]
\center
\begin{tabular}
{l  c c  c c  c c c c r} \hline \hline

RK1     & \multicolumn{2}{c}{RK2}       & RK3 & \multicolumn{2}{c}{RK4}   & $\ldots $ & RK7 & $\ldots $ & RKN    \\
Euler   &  midpoint    & trapezoid & Radau IIA & Hermite-Simpson  & classic & & &  & pseudospectral \\ [0.7ex]
\hline \\[-1.2ex]

 $ \mathcal{O}(h) $  & $ \mathcal{O}(h^2) $   &$ \mathcal{O}(h^2) $     & $ \mathcal{O}(h^3) $&
 $ \mathcal{O}(h^4) $  &  $ \mathcal{O}(h^4) $  & $\cdots$ & $ \mathcal{O}(h^7) $  & $\cdots$ & $ \mathcal{O}(h^N) $ \\[0.7ex]

 \hline \hline
\end{tabular}
\caption{Convergence rates from Euler to pseudospectral discretizations\cite{hnw-ode,fornberg,boyd,ap-ode}.}\label{table:convergence-rates}
\end{table*}
%
% ----------------------------------------------------------
In Table~\ref{table:convergence-rates} $h$ is the step size and the last column is written as $h^N$ for simplicity of presentation instead of $h^{N+1}$.  Table~\ref{table:convergence-rates} also implies that a PS discretization may be viewed as an extreme ``right'' end of a generalized Runge-Kutta method while an Euler method is at the left-most end of this spectrum (of finite bandwidth). Because there is nothing to the right of it, Table~\ref{table:convergence-rates} explains why a PS discretization (when properly implemented\cite{advances,DIDO:arXiv,ross-book,Radau-GNC05,auto-knots}) consistently outperforms other discretizations.

It is well-known\cite{hnw-ode,ap-ode,trefethen-2000,boyd} that the theoretical rates of convergence indicated in Table~\ref{table:convergence-rates} drops considerably for sufficiently small values of $h$ due to roundoff errors.  In Lagrange PS methods (i.e., those based on differentiation matrices) the onset of roundoff errors can be mitigated (but not avoided\cite{twist,boyd,trefethen-2000}).  A Birkhoff method does not suffer from this problem provided the grid points are computed accurately; hence, if one chooses a grid such that\cite{boyd,fornberg,newBirk-part-II,trefethen-2000},
\begin{equation}\label{eq:h=1overN}
h = \mathcal{O}(1/N)
\end{equation}
then, it follows that the theoretical convergence rate of a Birkhoff method can be obtained by substituting \eqref{eq:h=1overN} in the last column of Table~\ref{table:convergence-rates}.  Carrying out this exercise generates the theoretically achievable ``infinite order'' rate of convergence,
\begin{multline}\label{eq:inf-order}
\text{Birkhoff Convergence Rate (Theoretical)} =  \mathcal{O}(h^N) \\
= \mathcal{O}(1/N)^N = \mathcal{O}\left( N^{-N} \right)
\end{multline}
The practical convergence rate of a Birkhoff method is slower than \eqref{eq:inf-order} because of a multitude of factors.  Nonetheless, its convergence rate is still substantially faster relative to the non-PS discretization methods denoted in Table~\ref{table:convergence-rates}.  Part of the reason for the practical drop in the infinite-order convergence rate of a Birkhoff discretization is its interaction between Layers~1 and 2 of Fig.~\ref{fig:TrajOptPerspective} (and Layer~0 itself).

\subsection{Interactions Between Layers~1 and 2: Convergence and Accuracy}
Let $\by^* \in \real{N_y}$ be an exact solution to a continuous optimization problem. Barring extremely special cases, an optimization algorithm cannot produce $\by^*$ within a finite number of iterations. Instead, it can generate a number $\by^\sharp \in \real{N_y}$ such that,
\begin{equation}\label{eq:dely=eps}
\norm{\by^\sharp - \by^*} \le \epsilon > 0
\end{equation}
where $\norm{\cdot}$ is some appropriate norm that is consistent with the properties of the algorithm (Layer~2) and $\epsilon > 0$ is a given number (i.e., tolerance). If $\epsilon$ is chosen to be a number that is less than the achievable discretization error indicated in Table~\ref{table:convergence-rates}, then the resulting trajectory optimization method (i.e., Layers~1 and 2 of Fig.~\ref{fig:TrajOptPerspective} combined) will yield meaningless solutions\cite{cullum:1972,boris:AMP:2004,polak:book-1997} or false infeasibilities\cite{TAC:linearizable,ross:roadmap-2005,CMP:history}.  Furthermore, even if one were to ignore round-off errors, $\epsilon$ cannot be chosen to be equal to machine precision, $\epsilon_M$, even with a selection of very small values of $h$.  This is because the best case $\epsilon$ in \eqref{eq:dely=eps} is $\sqrt{\epsilon_M}$\cite{Gill:book,DIDO:arXiv,NW:NumOptBook,luenberger-2008}; hence, we must have,
\begin{equation}\label{eq:eps=sqrt}
\epsilon \ge \sqrt{\epsilon_M}
\end{equation}
To put this point in perspective, suppose $\epsilon_M = 10^{-16}$ (IEEE double-precision floating-point arithmetic) then \eqref{eq:eps=sqrt} requires $\epsilon \ge 10^{-8}$. Thus, any claims of accuracy greater than eight significant digits must be viewed with great suspicion. Regardless, a fast trajectory optimization method must effectively coordinate the convergence rates indicated in Table~\ref{table:convergence-rates} with the theoretical and practical achievable accuracies indicated in \eqref{eq:dely=eps} and  \eqref{eq:eps=sqrt}.  In the spectral algorithm\cite{spec-alg,auto-knots,DIDO:arXiv}, $\epsilon$ is set to $\epsilon^N$ (i.e., a grid-dependent number) such that,
\begin{equation}\label{eq:limeps=0}
\lim_{N \to \infty} \epsilon^N = 0
\end{equation}
Equation~\eqref{eq:limeps=0} is necessary for theoretical proofs of convergence\cite{cullum:1972,boris:AMP:2004,polak:book-1997} even though the minimum practical value of $\epsilon^N$ is bounded from below by \eqref{eq:eps=sqrt}.  To generate a faster spectral algorithm, it is necessary to moderate $\epsilon^N$ with the grid density ($\sim 1/N$) in a manner that is consistent with the convergence rate\cite{TAC:linearizable,kang-rate-2010,Kang_2008_convergence,PSReview-ARC-2012}.

Yet another point to consider for fast trajectory optimization is oracle complexity\cite{nesterov-book-2004,luenberger-2008} (i.e., number of function calls) to achieve a targeted accuracy ($\epsilon$ or $\epsilon^N$ as noted in the previous paragraph).  The oracle complexity (of Layer~2 in Fig.~\ref{fig:TrajOptPerspective}) depends upon the choice of the discretization method (Layer~1).  Furthermore, within a given oracle, there is the problem of computational complexity of returning the function value (and possibly its derivatives).  As shown in \cite{sandia-aas-23}, the computational complexity of returning a function value of a Birkhoff-discretized problem can be reduced from $\mathcal{O}(N^2)$ to $\mathcal{O}(N\log(N))$ using an FFT\cite{kreyszig-2011,fft-original,fftw}.  See also \cite{newBirk-part-II} for further discussions on oracle complexity.

There are far more issues to consider in generating a fast trajectory optimization method such as constructing a low-complexity method for solving the linear algebra that pairs with the chosen discretization method (Layer~1) and the algorithm of choice (Layer~2)\cite{sandia-aas-23}. In view of all such details that are presently unknown vis-a-vis a comprehensive theory of computational complexity for trajectory optimization, we use the word fast to denote the substantial reduction in computational complexity of the subcomponents of Fig.~\ref{fig:TrajOptPerspective}.

\section{Mathematical Foundations of Birkhoff Interpolants over an Arbitrary Grid}
\label{sec:NewBirkGeneric}

\begin{figure}[h!]
      \centering
      {\parbox{0.9\columnwidth}{
      \centering
      {\includegraphics[width = 0.95\columnwidth]{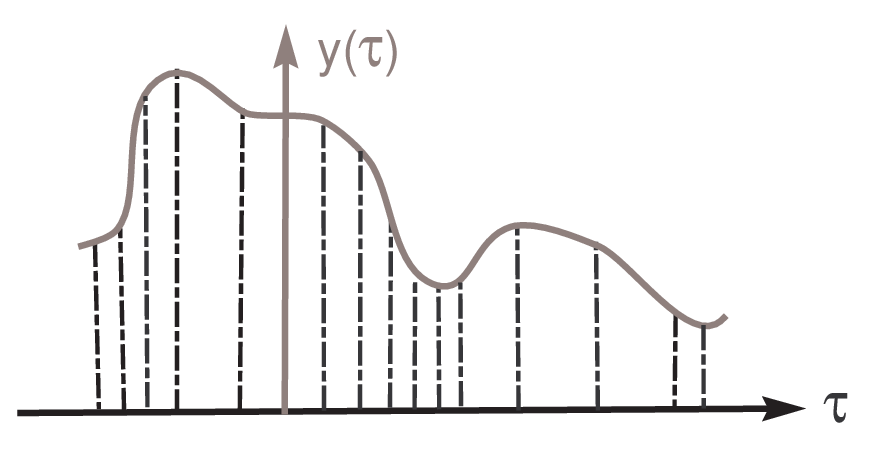}}
      \caption{\textsf{Schematic for an arbitrary grid. Figure adapted from \cite{furtherResults} with permission from I. M. Ross \copyright\ I. M. Ross, 2019.}}\label{fig:ArbGrid_y}
      }
      }
\end{figure}

As noted elsewhere\cite{fastmesh,furtherResults,arb-grid-aas, arb-grid}, the development of a PS method using a Gaussian grid at the outset tends to obfuscate an otherwise clear idea. Consequently, we follow \cite{furtherResults,arb-grid} and use an arbitrary grid as a starting point. More importantly, the use of an arbitrary grid and possibly non-orthogonal, non-polynomial basis functions\cite{nonpoly-interp,rbf-interp,sinc-interp,wavelets} generates a simpler and more elegant theory.  Furthermore, the use of an arbitrary grid for the development of a PS theory explains when and why a Gaussian condition or orthogonality is imposed\cite{arb-grid,furtherResults}. Note also that a Gaussian grid is neither necessary nor critical for convergence\cite{arb-grid}.  With these observations in focus, we begin by defining
\begin{equation}
\pi^N := \set{\tau_0, \tau_1, \ldots, \tau_N}, \quad \tau_k \in \Real, \ k \in \mathbb{N}
\end{equation}
to be a set of arbitrary grid points (see Fig.~\ref{fig:ArbGrid_y})  over the time interval $[\tau^a, \tau^b]$ such that
\begin{equation}\label{eq:arb-grid}
-\infty < \tau^a \le \tau_0 < \tau_1 < \cdots < \tau_{N-1} < \tau_N \le \tau^b < \infty
\end{equation}
%
%======================================================================================
%==========================================================================================
For the purposes of simplicity in the presentation of the ideas, we assume a finite horizon while noting that the results presented in \cite{Radau-GNC05,Radau-JGCD,LGR-rtoc} easily allow its extension to infinite-horizon problems; see also \cite{ross-book} and the references contained therein for additional details pertaining to solving infinite-horizon problems.

Let, $[\tau^a, \tau^b] \ni \tau \mapsto y(\tau) \in \Real$ be a given differentiable function. In \cite{newBirk-2023}, we defined two Birkhoff interpolants (over the same grid, $\pi^N$) given by,
\begin{subequations}\label{eq:Birk-interps}
\begin{align}
I^N_a y(\tau) &:=  y(\tau^a) B_0^0(\tau) + \sum_{j=0}^N \dot y(\tau_j) B_j^a(\tau)\label{eq:Birk-interp-a}\\
I^N_b y(\tau) &:= \sum_{j=0}^N \dot y(\tau_j) B_j^b(\tau) + y(\tau^b) B_N^N(\tau)\label{eq:Birk-interp-b}
\end{align}
\end{subequations}
where, $I^N_\theta, \theta \in \set{a, b}$ are the two interpolation operators and  $B_j^\theta, \ \theta \in \set{a, b}, j = 0,  \ldots, N$ are (Birkhoff) basis functions that satisfy the interpolation conditions,
\begin{align}\label{eq:birk-conditions}
\begin{aligned}
B_0^0(\tau^a)           &= 1,            &&    B_N^N(\tau^b)   = 1,                      \\
 {\dot B}_0^0(\tau_i)   &=  0,          &&   {\dot B}_N^N(\tau_i)   = 0,   & i = 0, \ldots, N \\
B_j^a(\tau^a)           &=  0,            && B_j^b(\tau^b) = 0,                       & j = 0, \ldots, N \\
{\dot B}_j^\theta(\tau_i) &=  \delta_{ij},  &&   i = 0, \ldots, N, & j = 0, \ldots, N
\end{aligned}
\end{align}
and $\delta_{ij}$ is the Kronecker delta.

%==================================================
\begin{remark}
Equation \eqref{eq:birk-conditions} is obtained by simply imposing the interpolation conditions,
\begin{subequations}
\begin{align}
I^N_a y(\tau^a) &= y(\tau^a),   &&\left.\frac{d}{d\tau} \left(I^N_a y(\tau)\right)\right|_{\tau=\tau_j} = \dot y(\tau_i),  &&i = 0, \ldots, N\\
I^N_b y(\tau^b) &= y(\tau^b),   &&\left.\frac{d}{d\tau} \left(I^N_b y(\tau)\right)\right|_{\tau=\tau_j} = \dot y(\tau_i),  &&i = 0, \ldots, N
\end{align}
\end{subequations}
\end{remark}
%-----------------------------
\begin{remark}\label{rem:no-poly-assumption}
There is no assumption of polynomials in all of the preceding equations. Hence, one may use radial basis functions\cite{rbf-interp}, Sinc functions\cite{sinc-interp}, wavelets\cite{wavelets} or any other basis function to generate Birkhoff interpolants.
\end{remark}
%-----------------------------
\begin{remark}
In principle, one can define an additional family of $(N+1)$ different Birkhoff interpolants over the exact same grid $\pi^N$ (and not just two) by replacing $y(\tau^a)$ in \eqref{eq:Birk-interp-a} by $y(\tau_k), k = 0, 1, \ldots, N$.
\end{remark}
%======================================
%
Unlike Lagrange or Hermite interpolants, a Birkhoff interpolant might not exist\cite{lorentz,schoenberg}; hence, it is critically important to prove its existence at the outset.
%====================
\begin{lemma}[Existence and Uniqueness Lemma\cite{newBirk-2023}]\label{prop:Birk-existence}
Let $L_j(\tau)$ be the antiderivative of $\ell_j(\tau)$, the $j$th Lagrange interpolating basis function over $\pi^N$.  Then, the Birkhoff basis functions that satisfy \eqref{eq:birk-conditions} are given explicitly by,
\begin{subequations}
\begin{align}
B^0_0(\tau) &= 1 =  B^N_N(x)                       &&     \label{eq:B00}\\
B_j^a(\tau) &= L_j(\tau) - L_j(\tau^a),                 && j = 0, \ldots, N \label{eq:Ba=Ldiff} \\
B_j^b(\tau) &= L_j(\tau) - L_j(\tau^b),                 && j = 0, \ldots, N \label{eq:Bb=Ldiff}
\end{align}
\end{subequations}
\end{lemma}
%================================
As noted in \cite{newBirk-2023}, these Birkhoff interpolants are different from the ones developed by Wang et al\cite{wang} but enjoy many of the same desirable properties.  They are shown to be universal in \cite{newBirk-2023} in the sense that they can be applied without change to solve high-order differential equations and a variety of boundary conditions, such as Dirichlet, Neumann and Robin.  It will be apparent shortly that they are also universal in their applications to trajectory optimization.
%=================================
\begin{definition}[\cite{newBirk-2023}]\label{def:birk-wts}
The Birkhoff quadrature weights, $w^B_j, \ j=0, \ldots, N$, are defined by,
\begin{equation}\label{eq:w-birk-def}
w^B_j := \int_{\tau_a}^{\tau_b}  \ell_j(\tau)\, d\tau
\end{equation}
\end{definition}
%----------------------------------
%===================================
\begin{lemma}[\cite{newBirk-2023}] \label{lemma:Ba-Bb=wts}
$\big(B^a_j(\tau) - B^b_j(\tau)\big)$ is a constant given by,
\begin{equation}\label{eq:Ba-Bb=wts}
B^a_j(\tau) - B^b_j(\tau) = w^B_j
\end{equation}
\end{lemma}
%======================================
%======================================
\begin{definition}\label{def:a=b}
The $a$- and $b$-forms of the Birkhoff interpolants for $y(\tau)$ given by \eqref{eq:Birk-interps} are said to be equivalent if $I^N_a y(\tau) = I^N_b y(\tau) \ \forall \tau \in [\tau^a, \tau^b]$.
\end{definition}
%======================================
\begin{proposition}\label{prop:a=b_iff_quad}
The $a$- and $b$-forms of the Birkhoff interpolants given by \eqref{eq:Birk-interps} are equivalent if and only if
\begin{equation}\label{eq:yN=y0+wtsV}
y(\tau^b) = y(\tau^a) + \sum_{j=0}^N w^B_j\, {\dot y}(\tau_j)
\end{equation}
\end{proposition}
%------------------------------------
\begin{proof}
We will first prove that if $I^N_a y(\tau) = I^N_b y(\tau) \ \forall \tau \in [\tau^a, \tau^b]$, then \eqref{eq:yN=y0+wtsV} holds.

Equating \eqref{eq:Birk-interp-a} and \eqref{eq:Birk-interp-b} we get,
\begin{equation}
 y(\tau^a) B_0^0(\tau) + \sum_{j=0}^N \dot y(\tau_j) B_j^a(\tau) =
\sum_{j=0}^N \dot y(\tau_j) B_j^b(\tau) + y(\tau^b) B_N^N(\tau)
\end{equation}
Equation~\eqref{eq:yN=y0+wtsV} then follows from Lemma~\ref{lemma:Ba-Bb=wts} and \eqref{eq:B00}.

To prove the converse, substitute \eqref{eq:yN=y0+wtsV} in \eqref{eq:Birk-interp-b}.  This generates,
\begin{align}
I^N_b y(\tau) &:= \sum_{j=0}^N \dot y(\tau_j) B_j^b(\tau) + y(\tau^b) B_N^N(\tau) \nonumber \\
&= \sum_{j=0}^N \dot y(\tau_j) B_j^b(\tau) + y(\tau^a)B_N^N(\tau) + \sum_{j=0}^N w^B_j\, {\dot y}(\tau_j) B_N^N(\tau) \nonumber \\
&= \sum_{j=0}^N \dot y(\tau_j) B_j^b(\tau) + y(\tau^a)B_0^0(\tau) + \sum_{j=0}^N (B^a_j(\tau) - B^b_j(\tau))\, {\dot y}(\tau_j) \label{eq:proof:a=b}
%
%& =y(\tau_0) B_0^0(\tau) + \sum_{j=0}^N \dot y(\tau_j) B_j^a(\tau) =
\end{align}
where, the last equation follows from Lemma~\ref{lemma:Ba-Bb=wts} and \eqref{eq:B00}.  Simplifying \eqref{eq:proof:a=b} yields $I^N_b y(\tau) = y(\tau^a) B_0^0(\tau) + \sum_{j=0}^N B^a_j(\tau) \, \dot y(\tau_j) = I^N_a y(\tau)$, by definition (cf.~\eqref{eq:Birk-interp-a}).
\end{proof}
%====================================
\begin{definition}[\cite{newBirk-2023}]\label{def:Birk-matrix}
Let $\theta \in \set{a, b}$.  The Birkhoff matrix, $\bB^\theta$, is defined by,
\begin{equation}\label{eq:Birk-theta-def}
\bB^\theta =
\begin{bmatrix}
B_0^\theta(\tau_0) & B_1^\theta(\tau_0) & \dots & B_N^\theta(\tau_0) \\
B_0^\theta(\tau_1) & B_1^\theta(\tau_1) & \dots & B_N^\theta(\tau_1) \\
\vdots & \vdots & &\vdots \\
B_0^\theta(\tau_N) & B_1^\theta(\tau_N) & \dots & B_N^\theta(\tau_N)
\end{bmatrix}
\end{equation}
\end{definition}
%----------------------------------
%
\begin{definition}\label{def:WB}
\begin{enumerate}\item[]
\item The $(N+1)$ Birkhoff quadrature weight vector is defined by  $\bw_B := [w^B_0, \ldots, w^B_N]^T$
\item The $(N+1) \times (N+1)$ Birkhoff quadrature weight matrix is defined by  $W_B := \textrm{diag}(w^B_0, \ldots, w^B_N)$
\end{enumerate}
\end{definition}
%==================================
\begin{lemma}[\cite{newBirk-2023} ]\label{lemma:lastRow=wts}
Suppose $\pi^N$ is such that $\tau_0 = \tau^a$ and $\tau_N = \tau_b$.  Then,
\begin{enumerate}
%\item[]
\item The last row of $\bB^a$ is identically equal to the Birkhoff quadrature weights.
\item The first row of $\bB^b$ is identically equal to the negative of the Birkhoff quadrature weights.
\end{enumerate}
\end{lemma}
%----------------------------------
%==================================
\begin{remark}
Lemma~\ref{lemma:lastRow=wts} implies that if $\pi^N$ is selected such that $\tau_0 = \tau^a$ and $\tau_N = \tau_b$, then the Birkhoff quadrature weights can be simply extracted out as the last row of $\bB^a$ or the negative of the first row of $\bB^b$. That is, no separate quadrature computations are necessary.
\end{remark}
%==================================
\begin{lemma}\label{lemma:pre-int-by-parts}
Let $\epsilon^N_{jk}$ be defined by
\begin{align}
\epsilon^N_{jk}  &:= \sum_{i=0}^N L_j(\tau_i) \ell_k(\tau_i) w^B_i - \int_{\tau^a}^{\tau^b} L_j(\tau)\ell_k(\tau)\, d\tau \label{eq:quad-error}
%\\ & =  L_j(\tau_k) w^B_k - \int_{\tau_0}^{\tau_N} L_j(\tau)\ell_k(\tau)\, d\tau
\end{align}
Then,
\begin{equation}\label{eq:pre-int-by-parts}
W_B \bB^b + [\bB^a]^T W_B = \bdelta^N
\end{equation}
where $\bdelta^N$ is an $(N+1) \times (N+1)$ symmetric matrix whose constituent elements, $\delta^N_{jk}$, are given by $\delta_{jk}^N = \epsilon^N_{jk} + \epsilon^N_{kj}$.
\end{lemma}
%---------------------------------
\begin{proof}
By chain rule, we have,
\begin{equation}\label{eq:chain-rule-inProof}
\frac{d}{d\tau} \Big(L_k(\tau) L_j(\tau)  \Big) = L_k(\tau) \ell_j(\tau) + \ell_k(\tau) L_j(\tau)
\end{equation}
Integrating both sides of \eqref{eq:chain-rule-inProof} over the interval $[\tau^a, \tau^b]$, we get,
\begin{equation}\label{eq:proof-int-by-parts}
L_k(\tau^b) L_j(\tau^b) - L_k(\tau^a) L_j(\tau^a)  = \int_{\tau^a}^{\tau^b} \Big( L_k(\tau) \ell_j(\tau) + \ell_k(\tau) L_j(\tau) \Big)\, d\tau
\end{equation}
From Definition~\ref{def:birk-wts} it follows that $L_j(\tau^b) = L_j(\tau^a) + w^B_j $ for all $j= 0, 1, \ldots, N$.  Substituting this relationship in the left hand side of \eqref{eq:proof-int-by-parts}, we get,
\begin{multline}
L_k(\tau^b) L_j(\tau^b) - L_k(\tau^a) L_j(\tau^a)  \\
 =  \big(L_k(\tau^a) + w^B_k \big) \Big( L_j(\tau^a) + w^B_j \Big) - L_k(\tau^a) L_j(\tau^a)  \\
 = L_k(\tau^a) w^B_j + L_j(\tau^a) w^B_k + w^B_k w^B_j  \\
 = L_k(\tau^a) w^B_j + L_j(\tau^b) w^B_k \label{eq:proof-lhs=}
\end{multline}
By hypothesis, the right-hand-side of \eqref{eq:proof-int-by-parts} can be written as,
\begin{multline}
\int_{\tau^a}^{\tau^b} \Big( L_k(\tau) \ell_j(\tau) + \ell_k(\tau) L_j(\tau) \Big)\, d\tau \\
= \sum_{i=0}^N  \Big( L_k(\tau_i) \ell_j(\tau_i) + \ell_k(\tau_i) L_j(\tau_i) \Big) w^B_i  - \epsilon^N_{kj} - \epsilon^N_{jk} \\
 = L_k(\tau_j) w^B_j + L_j(\tau_k) w^B_k - \delta_{kj}^N \label{eq:proof-rhs=}
\end{multline}
Equating the right-hand-sides of \eqref{eq:proof-lhs=} and \eqref{eq:proof-rhs=} we get,
\begin{equation}\label{eq:proof-pre-result}
\Big( L_k(\tau_j) -  L_k(\tau^a) \Big) w^B_j + \Big( L_j(\tau_k) - L_j(\tau^b) \Big) w^B_k - \delta^N_{kj} = 0
\end{equation}
From Lemma~\ref{prop:Birk-existence}, it follows that \eqref{eq:proof-pre-result} can be rewritten as,
\begin{equation}\label{eq:proof-result=almost}
 B^a_k(\tau_j) w^B_j + B^b_j(\tau_k) w^B_k = \delta^N_{kj}
\end{equation}
Using Definitions~\ref{def:Birk-matrix} and \ref{def:WB}, it follows that the left-hand-side of \eqref{eq:proof-result=almost} is the $kj$-th element of \eqref{eq:pre-int-by-parts} for $k, j = 0, 1, \ldots, N$.
\end{proof}
%===================================

\section{Birkhoff Discretizations of Optimal Control Problems}
\label{sec:direct=PS+S}

In the development of spectral and pseudospectral methods for optimal control\cite{arb-grid}, it is frequently sufficient to generate the foundational equations and theorems for a ``distilled problem'' that captures all of the key elements necessary for generalization; see for example, \cite{fastmesh,furtherResults,arb-grid,advances}.  Once this is done, the production of the requisite computational equations for a general optimal control problem becomes straightforward albeit with a laborious exercise in bookkeeping\cite{furtherResults,knots,hybrid:jgcd,acc:hybrid}. Management techniques for simplifying the bookkeeping are provided in \cite{newBirk-part-II} with many details described elsewhere\cite{furtherResults,hybrid:jgcd,DIDO:arXiv}. Because the clutter of bookkeeping obscures the main ideas in producing the foundational equations for a Birkhoff theory for trajectory optimization, we follow \cite{advances} and consider the following distilled nonlinear optimal control problem:
%
%===============================================
\begin{eqnarray}
& (P) \left\{
\begin{array}{lrl}
\emph{Minimize } && J[x(\cdot), u(\cdot)] := E(x(\tau^a), x(\tau^b))\\[0.5em]
\emph{Subject to}&& \dot x(\tau) = f(x(\tau), u(\tau)) \\
                 && e(x(\tau^a), x(\tau^b)) \le  0
\end{array} \right. & \label{eq:probP}
\end{eqnarray}
%===============================================
In \eqref{eq:probP},
$J$ is a cost functional given by a differentiable endpoint cost function $E: \Real \times \Real \to \Real$, the pair $\big(x(\cdot), u(\cdot)\big)$ is the unknown system trajectory (i.e., state-control function pair, $\tau \mapsto \big(x(t), u(t)\big) \in \Real \times \Real$), $f: \Real \times \Real \to \Real$ is a given differentiable dynamics function, and $e: \Real \times \Real \to \Real$ is a given differentiable endpoint constraint function that constrains the endpoint values of $x(\tau)$.

As a starting point, consider first the use \eqref{eq:Birk-interp-a} to represent the state trajectory over an arbitrary grid $\pi^N$. Let,
\begin{equation}\label{eq:xN=Birk-a-1D}
x^N(\tau) := x^a B_0^0(\tau) + \sum_{j=0}^N v_j B_j^a(\tau)
\end{equation}
where, $x^a:= x^N(\tau^a)$ and $v_j:= d/d\tau (x^N(\tau))_{\tau_j}, \ j = 0, 1, \ldots, N$ are all unknowns.
%---------------------------------------
\begin{definition}\label{def:virtual-variables}
Following \cite{DIDO:arXiv}, we call $v_j, \ j = 0, 1, \ldots, N$ the virtual variables.
\end{definition}
%---------------------------------------
\begin{remark}\label{rem:xN=semi-analytic}
If $x^a$ and $v_j,\ j = 0, 1, \ldots, N$ are known, \eqref{eq:xN=Birk-a-1D} represents a closed-form expression of a state trajectory in terms of the Birkhoff basis functions.  Consequently, we may consider the Birkhoff theory as a means to produce semi-explicit or semi-discrete solutions to an optimal control problem.
\end{remark}
Substituting \eqref{eq:xN=Birk-a-1D} in the dynamical equation, $\dot x(\tau) = f(x(\tau), u(\tau))$, we can define a residual function, $R^N_{\dot x}$, according to:
\begin{equation}\label{eq:residual-xdot}
R^N_{\dot x}(x^N(\tau), u(\tau)):= \dot x^N(\tau) - f(x^N(\tau), u(\tau))
\end{equation}
At first blush, setting $R^N_{\dot x}(x^N(\tau), u(\tau)) = 0$ for all $\tau \in [\tau^a, \tau^b]$ seems highly desirable; however, as is well-recognized (see \cite{cullum:1972,TAC:linearizable,halkin:1966,boris:AMP:2004,CMP:history}) this is neither warranted nor preferred either from a mathematical or an engineering point of view.  In constructing a viable computational optimal control theory, we seek to find state-control function pairs, $\big( x^N(\cdot), u(\cdot) \big)$, where $R^N_{\dot x}(x^N(\tau), u(\tau))$ vanishes when projected along some preferred basis ``test'' functions\cite{boyd}, $\phi_i(\tau), \ i = 0, \ldots, N$.  That is, we set
\begin{multline}\label{eq:xdot-f=0viaMWR}
\langle R^N_{\dot x}, \phi_i \rangle_\varrho := \int_{\tau^a}^{\tau^b} R^N_{\dot x}\left(x^N(\tau), u(t) \right) \phi_i(\tau)\varrho(\tau)\, d\tau = 0, \\
i = 0, \ldots, N
\end{multline}
where, $\langle \cdot, \cdot \rangle_\varrho$ denotes an inner product with  $\varrho$ as the weight function. Specific choices of $\phi_i(\tau)$ (and $\varrho(\tau)$) generate specific projections. To limit the scope of this paper, we set $\varrho(\tau) = 1$, a constant.  Furthermore, we will eventually limit the choice of $\phi_i(\tau)$ to two classes of functions corresponding to a new unified theory for spectral and pseudospectral methods for trajectory optimization.

\subsection{Development of Birkhoff Pseudospectral Discretizations of Problem~$(P)$}\label{subsec:PNa+b}

Taking $\phi_i(\tau) = \delta^D(\tau - \tau_i)$, where $\delta^D(\tau - \tau_i)$ are Dirac delta functions centered on the grid points, $\tau_i, \ i = 0, \ldots, N$, \eqref{eq:xdot-f=0viaMWR} generates,
\begin{equation}\label{eq:vj=fj-1D}
\langle R^N_{\dot x}, \delta^D(\tau - \tau_i) \rangle_{\varrho=1} = 0  \Rightarrow  v_i = f(x^N(\tau_i), u(\tau_i)), \quad i = 0, \ldots, N
\end{equation}
It is apparent that \eqref{eq:vj=fj-1D} is identical to imposing the state derivative condition over all grid points.  This is the pseudospectral version of \eqref{eq:xdot-f=0viaMWR}. That is, \eqref{eq:xdot-f=0viaMWR} is the foundational equation, and \eqref{eq:vj=fj-1D} was derived from this general principle.
%--------------------
\begin{remark}
The controls, $u(\tau_i), \ i = 0, 1, \ldots, N$, that solve \eqref{eq:vj=fj-1D} are not necessarily of polynomial origin even if the Birkhoff basis functions, $B^a_j(\tau)$ are constructed using polynomials.  This point follows from the fact that \eqref{eq:vj=fj-1D} was obtained without ever assuming $t \mapsto u(t)$ was given by an equivalent version of \eqref{eq:xN=Birk-a-1D} or some other interpolating polynomial. As a result, \eqref{eq:vj=fj-1D} may also be viewed as part of a defining equation for $u(\tau)$ given implicitly by the inverse map (``$ f^{-1}$'') of the dynamics function, $f$. Because $f^{-1}$ is frequently not a polynomial, neither is $u(\cdot)$. Hence, $u(\tau_i), \ i = 0, 1, \ldots, N$ are fully capable of representing discontinuities without suffering any Gibbs phenomenon.  This point is illustrated in \cite{newBirk-part-II}.
\end{remark}
%--------------------
Analogous to \eqref{eq:residual-xdot}, define,
\begin{equation}
R^N_{x}(x^N(\tau), u(\tau)):=  x^N(\tau)- x^a B_0^0(\tau) + \sum_{j=0}^N v_j B_j^a(\tau)
\end{equation}
Setting $\langle R^N_{x}, \delta^D(\tau - \tau_i) \rangle_{\varrho=1} = 0, \ i = 0, \ldots, N$ generates,
\begin{equation}\label{eq:xi=xa+Bijvj-1d}
x^N(\tau_i) = x^a  + \sum_{j=0}^N v_j B_j^a(\tau_i), \quad i = 0, \ldots, N
\end{equation}
Let,
\begin{multline}\label{eq:XVandUbydef-1D}
X := [x^N(\tau_0), \ldots, x^N(\tau_N)]^T, \quad  V := [v_0, \ldots, v_N]^T,  \\
U := [u(\tau_0), \ldots, u(\tau_N)]^T
\end{multline}
Then, \eqref{eq:vj=fj-1D} and \eqref{eq:xi=xa+Bijvj-1d} can be vectorized as,
\begin{subequations}\label{eq:vectorizedXVeqn}
\begin{align}
X &= x^a \bb + \bB^a V \label{eq:X=Birk-a-1D}\\
V &= f(X, U) \label{eq:V=f-vec-1D}
\end{align}
\end{subequations}
where, $\bb$ is an $(N+1)\times 1$ vector of ones and $f(X, U) := [f(x^N(\tau_0), u(\tau_0)), \ldots, f(x^N(\tau_N), u(\tau_N))]^T $.
\begin{remark}\label{rem:0=eps}
\begin{enumerate}\item[]
\item Strictly speaking, \eqref{eq:vectorizedXVeqn} must be written more technically in terms of an additive $\epsilon^N$ relaxation with $\epsilon^N \to 0$ as $N \to \infty$; see \eqref{eq:limeps=0}.  This is because, it has long been known\cite{ross:roadmap-2005,polak:book-1997,TAC:linearizable,cullum:1972,boris:AMP:2004} that if the discretized equations of an optimal control problem are imposed exactly (as implied in \eqref{eq:vectorizedXVeqn}), then the resulting set of equations may not even produce a feasible solution even if the original problem had infinitely many feasible solutions. Simple (counter) examples are provided in \cite{CMP:history,boris:AMP:2004,cullum:1972,halkin:1966} to illustrate this point.  To facilitate a simpler presentation, we implicitly assume that the equality in \eqref{eq:vectorizedXVeqn} is relaxed to conform with well-established theoretical results in computational optimal control.
 \item In practical trajectory optimization, \eqref{eq:vectorizedXVeqn} is imposed approximately by way of primal feasibility tolerances\cite{arb-grid,DIDO:arXiv,spec-alg}.  Furthermore, as noted in the discussions pertaining to \eqref{eq:eps=sqrt}, these tolerances are required to be no smaller than $\sqrt{\epsilon_M}$\cite{Gill:book,DIDO:arXiv}.

\end{enumerate}
\end{remark}
%----------------------------
In view of Remark~\ref{rem:0=eps}, we caveat all equality constraints in both parts I and II of this paper by an implicit assumption of relaxation.  Keeping this technicality in the background allows us to maintain our exposition accessible to a broader audience; nonetheless,
because this point is quite important, we memorialize this historical knowledge base in terms of the following definition that we implicitly use throughout this paper:
%--------------
\begin{definition}\label{def:approx=exact}
We say $x^N = y^N$ for a given pair  $\set{x^N, y^N}, N \in \mathbb{N}$ if, for any given $\epsilon > 0$, there exists an $N_\epsilon \in \mathbb{N}$ such that for all $N \ge N_\epsilon$, $\abs{x^N-y^N} \le \epsilon$.
\end{definition}
%----------------------

To complete a Birkhoff PS representation of Problem~$(P)$, we now make the first of two hypotheses in the selection of a grid $\pi^N$:
%
%===========================
\begin{hypothesis}\label{hyp:a=b}
There exists a grid $\pi^N$ and an $N_\epsilon \in \mathbb{N}$ such that for all $N \ge N_\epsilon$ the $a$- and $b$-forms of the Birkhoff interpolants are equivalent (in the sense of Definition~\ref{def:approx=exact}).
\end{hypothesis}
%=============================
%
Using Hypothesis~\ref{hyp:a=b}, Proposition~\ref{prop:a=b_iff_quad} and \eqref{eq:xN=Birk-a-1D}, we now employ the Birkhoff grid-equivalency condition,
\begin{equation}\label{eq:xb=xa+w*V}
x^b = x^a + \bw_B^T V
\end{equation}
to impose boundary conditions on $x^N(\tau^b) := x^b $.  Collecting all the relevant equations, we define Problem~$(P^N_a)$ as follows:
%
%===========================================================
\begin{align}
X \in \real{N+1}, \quad U \in \real{N+1}, \quad V \in \real{N+1}, \quad (x^a, x^b) \in \real{2}   & \nonumber\\
 \left(\textsf{$P_a^N$}\right) \left\{
\begin{array}{lrl}
\emph{Minimize } & J^N_a[X, U, V, x^a, x^b] :=& E(x^a,x^b)\\
\emph{Subject to} & X  = &x^a\,\bb + \bB^a V \\
& V =& f(X, U)  \\
& x^b = & x^a + \bw_B^T V \\
& e(x^a, x^b)  \le & 0
\end{array} \right. & \label{eq:Prob-PNa}
\end{align}
%===============================================
\begin{remark}\label{rem:lin-nonlin-split}
Structurally, Problem~$(P^N_a)$ comprises two components:
%----------------
\begin{enumerate}
\item A problem-invariant linear (and hence, convex) system of equations that depends only on the choice of $\pi^N$, and
\item A problem-dependent system whose evaluations require a computation of the data functions (i.e., $E, f$ and $e$) over  $\pi^N$.
\end{enumerate}
%----------------
\end{remark}
%--------------------------------------------------
\begin{remark}
As defined here in \eqref{eq:Prob-PNa}, Problem~$(P^N_a)$ is structurally similar to the ones defined in \cite{fastmesh} and \cite{furtherResults}.  The major difference between these problem definitions is in the definition of the Birkhoff matrix itself. The other difference is in the inclusion of the Birkhoff grid-equivalency condition given by \eqref{eq:xb=xa+w*V}.
\end{remark}
%----------------------------------------------

Using \eqref{eq:Birk-interp-b} as a starting point (instead of \eqref{eq:Birk-interp-a}) and following the same process as in the previous paragraphs, it is straightforward to show that an equivalent Problem~$(P^N_b)$ can be generated according to,
%
%===========================================================
\begin{align}
X \in \real{N+1}, \quad U \in \real{N+1}, \quad V \in \real{N+1}, \quad (x^a, x^b) \in \real{2}   & \nonumber\\
 \left(\textsf{$P_b^N$}\right) \left\{
\begin{array}{lrl}
\emph{Minimize } & J^N_b[X, U, V, x^a, x^b] :=& E(x^a,x^b)\\
\emph{Subject to} & X  = &x^b\,\bb + \bB^b V \\
& V =& f(X, U)  \\
& x^b = & x^a + \bw_B^T V \\
& e(x^a, x^b)  \le & 0
\end{array} \right. & \label{eq:Prob-PNb}
\end{align}
%===============================================
\begin{remark}
If $\pi^N$ is chosen according to Lemma~\ref{lemma:lastRow=wts}, then $\bw_B$ can be simply extracted out from  $\bB^\theta,\ \theta \in \set{a, b}$ and no additional computations are necessary to construct the linear Birkhoff system (cf.~Remark~\ref{rem:lin-nonlin-split}).
\end{remark}
%==========================

\subsection{Development of Birkhoff Spectral Approximations to Problem~$(P)$}\label{subsec:PNa+b-stars}
To produce a spectral-Galerkin\cite{boyd}  representation of Problem~$(P)$, consider the projection equation $\langle R^N_{\dot x}, \phi_i(\tau) \rangle_{\varrho=1} = 0$ given by,
\begin{multline}\label{eq:xdot=f-Gal-def}
\langle R^N_{\dot x}, \phi_i(\tau) \rangle_{\varrho=1} := \int_{\tau^a}^{\tau^b}\left(\dot x^N(\tau) - f(x^N(\tau), u(\tau)) \right) \phi_i(\tau)\, d\tau = 0, \\
i = 0, \ldots, N
\end{multline}
Substituting the Birkhoff representation of the state trajectory (i.e., \eqref{eq:xN=Birk-a-1D}) in \eqref{eq:xdot=f-Gal-def} we get,
\begin{equation}\label{eq:xdot=f-Gal-N}
\int_{\tau^a}^{\tau^b}\left(\sum_{j=0}^N v_j \ell_j(\tau) - f\big(x^N(\tau), u(\tau)\big) \right) \phi_i(\tau)\, d\tau = 0, \quad i = 0, \ldots, N
\end{equation}
To evaluate the left-hand-side of \eqref{eq:xdot=f-Gal-N} requires problem-specific information (i.e., $f(x,u)$).  Even if this is available, it may not be feasible to produce closed-form expressions of the required integrals. In order to devise a general purpose scheme, we now make the second (of two) hypotheses in choosing $\pi^N$:
%
%============================
\begin{hypothesis}\label{hyp:piN=quadPts}
Let $\tau \mapsto y(\tau) \in \Real$ be a bounded, integrable function over $[\tau^a, \tau^b]$.  Let,
\begin{equation}
Q^N(y(\cdot)) := \sum_{i=0}^N y(\tau_i) \, w^B_i - \int_{\tau^a}^{\tau^b} y(\tau)\, d\tau
\end{equation}
Then, given any $\epsilon > 0$, there exists a $\pi^N$ and an $N_\epsilon$ such that $\abs{Q^N(y(\cdot))} \le \epsilon$ for all $N \ge N_\epsilon$.
\end{hypothesis}
%============================
Using Hypothesis~\ref{hyp:piN=quadPts}, \eqref{eq:xdot=f-Gal-N} can be written as,
\begin{multline}\label{eq:xdot=f-Gal-wQuad}
\sum_{k=0}^N w^B_k \left(\sum_{j=0}^N v_j \ell_j(\tau_k) - f(x^N(\tau_k), u(\tau_k)) \right) \phi_i(\tau_k) = 0, \\
 i = 0, \ldots, N
\end{multline}
Because the quadrature in \eqref{eq:xdot=f-Gal-wQuad} approximates the left-hand-side of \eqref{eq:xdot=f-Gal-N}, the equality in the former equation must be viewed in the context of Definition~\ref{def:approx=exact}. As noted in Remark~\ref{rem:0=eps} we implicitly assume such relaxations. Note, in particular, that \eqref{eq:xdot=f-Gal-wQuad} offers a natural relaxation of its original equation (namely, \eqref{eq:xdot=f-Gal-N}) which is surprisingly more helpful (cf.~Remark~\ref{rem:0=eps}) than its exact imposition\cite{ross:roadmap-2005}.

For a spectral representation of the equations that constitute Problem~$(P)$, the basis test functions $\phi_i(\tau), \ i = 0, \ldots, N$ are chosen to be the same as $\ell_i(\tau), \ i = 0, \ldots, N$.  As a result, \eqref{eq:xdot=f-Gal-wQuad} simplifies dramatically to:
\begin{equation}\label{eq:xdot=f-Gal-approx}
w^B_i v_i  - w^B_i  f(x^N(\tau_i), u(\tau_i))   = 0, \quad i = 0, \ldots, N
\end{equation}
Vectorizing \eqref{eq:xdot=f-Gal-approx}, the Birkhoff spectral approximation to the dynamical equation can be written compactly as,
\begin{equation}\label{eq:xdot=f-Gal-hadamard}
\bw_B \circ V - \bw_B \circ f(X, U)  = 0
\end{equation}
where, $\circ$ represents the Hadamard product.

In performing a similar exercise in discretizing \eqref{eq:xN=Birk-a-1D} via the Galerkin method, it can be readily verified that we get the following spectral representation for state parameters:
\begin{equation}
x^a \bw_B \circ \bb + \bw_B \circ (\bB^a V)  - \bw_B \circ X = 0
\end{equation}
Collecting all the relevant equations, the Birkhoff spectral discretization of Problem~$(P)$ can be defined according to the following:
%===========================================================
\begin{align}
X \in \real{N+1}, \quad U \in \real{N+1}, \quad V \in \real{N+1}, \quad (x^a, x^b) \in \real{2}   & \nonumber\\
 \left(\textsf{$P_{a^*}^N$}\right) \left\{
\begin{array}{lll}
\emph{Minimize } & J^N_a[X, U, V, x^a, x^b] := E(x^a,x^b)\\
\emph{Subject to} & x^a \bw_B \circ \bb + \bw_B \circ \left(\bB^a V\right)  - \bw_B \circ X =  0 \\
& \bw_B \circ V - \bw_B \circ f(X, U)  =  0  \\
& x^b - x^a - \bw_B^T V =  0 \\
& e(x^a, x^b)  \le  0
\end{array} \right. & \label{eq:Prob-PNa*}
\end{align}
%===============================================
\begin{remark}
Comparing \eqref{eq:Prob-PNa*} to \eqref{eq:Prob-PNa}, it is apparent that Problem~$(P^N_{a^*})$ is very similar to Problem~$(P^N_{a})$.  The major difference between these problem definitions is the appearance of
Hadamard products with the Birkhoff quadrature weights in \eqref{eq:Prob-PNa*}.  The implication of the appearance of these quadrature weights will be apparent in Section~\ref{sec:CMP-derive}.
\end{remark}
%-----------------------------------------
\begin{remark}\label{rem:fourPNs}
Analogous to \eqref{eq:Prob-PNb}, one can similarly define Problem~$(P^N_{b^*})$.  These details are not provided for the purposes of brevity. It thus follows that Problem~$(P)$ may be discretized in at least four different ways denoted by Problems~$(P^N_{a}), (P^N_{b}), (P^N_{a^*})$ and $(P^N_{b^*})$.
\end{remark}
%=========================================

\section{Birkhoff Discretizations of Differential-Algebraic BVPs with Complementarity Conditions} \label{sec:16PlamNs}
Applying Pontryagin's Principle\cite{ross-book}, it is straightforward to show that the collection of necessary conditions for Problem~$(P)$ can be summarized in terms of the following differential-algebraic BVP with complementarity boundary conditions:
\begin{eqnarray}
& (\textsf{$P^\lambda$}) \left\{
\begin{array}{lrl}
& \dot x(\tau) =& f(x(\tau), u(\tau))  \\
&-\dot\lambda(\tau) = & \partial_x f(x(\tau), u(\tau)) \lambda(\tau)\\
& 0 = &\partial_u f(x(\tau), u(\tau)) \lambda(\tau) \\
& 0 \ge & e(x(\tau^a), x(\tau^b))  \\
& -\lambda(\tau^a) = & \partial_{x^a} \overline{E}(\nu, x^a, x^b) \\
& \lambda(\tau^b) = & \partial_{x^b} \overline{E}(\nu, x^a, x^b) \\
&  \nu\ \ \dagger  &e(x(\tau^a), x(\tau^b))
\end{array} \right. & \label{eq:probPlam}
\end{eqnarray}
where, $\lambda(\tau) \in \Real $ is the costate, $\overline{E}(\nu, x^a, x^b)$  is the endpoint Lagrangian function defined by,
\begin{equation}
\overline{E}(\nu, x^a, x^b) :=  E(x^a, x^b) + \nu\, e(x^a, x^b)
\end{equation}
and $\nu \in \Real$ is an endpoint multiplier that satisfies the complementarity condition abbreviated as $\nu\ \dagger  e(x(\tau^a), x(\tau^b))$\cite{ross-book}.

%-------------
\begin{remark}
It is apparent that Problem~$(P^\lambda)$ is well beyond a classic ``two-point'' BVP\cite{brysonHo}; in fact, it contains many of the acknowledged difficulties\cite{longuski,conway:survey,brysonHo} in solving trajectory optimization problems such as differential-algebraic constraints and complementarity conditions at the boundary points.  Hence, calling Problem~$(P^\lambda)$ simply as a BVP (or a two-point BVP) is quite misleading in the sense that it diminishes the combinatorial complexity of the complementarity conditions.

\end{remark}
%-------------

Suppose we use \eqref{eq:Birk-interp-b} to represent the costate trajectory over the arbitrary grid $\pi^N$.  In following the same procedure as in Section~\ref{sec:direct=PS+S}, we define,
\begin{equation}\label{eq:lamN=Birk-b-1D}
\lambda^N(\tau) := \lambda^b B_N^N(\tau) + \sum_{j=0}^N \varpi_j B_j^b(\tau)
\end{equation}
where, $\lambda^b := \lambda^N(\tau^b)$ and $\varpi_j:= d/d\tau (\lambda^N(\tau))_{\tau_j}, \ j = 0, 1, \ldots, N$ are all unknowns. Similar to the observation of Remark~\ref{rem:xN=semi-analytic}, we note the following:
%-----------------------
\begin{remark}
If $\lambda^b$ and $\varpi_j\ j = 0, 1, \ldots, N$ are known, \eqref{eq:lamN=Birk-b-1D} represents a closed-form expression of a costate trajectory.  As a result, a Birkhoff theory for solving Problem~$(P^\lambda)$ produces semi-explicit or semi-discrete solutions to the state-costate trajectory pair.
\end{remark}
%------------------------
%---------------------------------------
\begin{definition}\label{def:covirtual-variables}
Similar to Definition~\ref{def:virtual-variables}, we follow \cite{DIDO:arXiv} and call $\varpi_j, \ j = 0, 1, \ldots, N$ the co-virtual variables.
\end{definition}
%---------------------------------------

Following the same process as in Section~\ref{sec:direct=PS+S}, we begin by developing the projection equation $\langle R^N_{\dot \lambda}, \phi_i(\tau) \rangle_{\varrho=1}$ and setting it to 0; this yields,
\begin{multline}\label{eq:adjoint-Gal-def}
\langle R^N_{\dot \lambda}, \phi_i(\tau) \rangle_{\varrho=1} := \\
\int_{\tau^a}^{\tau^b}\left(\dot\lambda^N(\tau) + \partial_x f(x^N(\tau), u(\tau)) \lambda^N(\tau)\right) \phi_i(\tau)\, d\tau = 0, \\
i = 0, \ldots, N
\end{multline}
Selecting the test function $\phi_i(\tau)$ to be the same as $\ell_i(\tau)$ and using Hypothesis~\ref{hyp:piN=quadPts} it can be easily verified that \eqref{eq:adjoint-Gal-def} simplifies analogous to \eqref{eq:xdot=f-Gal-approx} to,
\begin{equation}\label{eq:adjoint-Gal-approx}
w^B_i \varpi_i  + w^B_i \partial_x f(x^N(\tau_i), u(\tau_i)) \lambda^N(\tau_i)  = 0, \quad i = 0, \ldots, N
\end{equation}
Let,
\begin{equation}\label{eq:LamOmbydef-1D}
\Lambda:= [\lambda^N(\tau_0), \ldots, \lambda^N(\tau_N)]^T, \quad  \Omega := [\varpi_0, \ldots, \varpi_N]^T
\end{equation}
Then, the spectral-Galerkin-Birkhoff approximation to the adjoint equation can be written compactly as,
\begin{equation}
\bw_B \circ \Omega + \bw_B \circ \, \partial_x f(X, U) \circ \Lambda = 0
\end{equation}
where, $\partial_x f$ is reused as an overloaded operator (see \cite{ross-book}, page 8) defined by,
$$ \partial_x f(X,U) := \big[\partial_x f(x^N(\tau_0), u(\tau_0)), \ldots, \partial_x f(x^N(\tau_N), u(\tau_N))\big]^T$$
In performing a similar exercise on imposing \eqref{eq:lamN=Birk-b-1D} via the spectral-Galerkin method, it can be readily shown that we get the following spectral-Galerkin-Birkhoff discretization of the costate:
\begin{equation}
\lambda^b \bw_B \circ \bb + \bw_B \circ (\bB^b \Omega)  - \bw_B \circ \Lambda = 0
\end{equation}
Collecting all the relevant equations, we now define a spectral-Birkhoff representation of Problem~$(P^\lambda)$ (a differential-algebraic BVP with complementarity boundary conditions) in terms of the following generalized root-finding problem:
\begin{eqnarray}
& \left(\textsf{$P_{a^*,b^*}^{\lambda, N}$}\right) \left\{
\begin{array}{lrl}
 & x^a \bw_B \circ \bb + \bw_B \circ \left(\bB^a V\right)  - \bw_B \circ X = & 0 \\
&\lambda^b \bw_B \circ \bb + \bw_B \circ \left(\bB^b \Omega \right)  - \bw_B \circ \Lambda= & 0  \\
& \bw_B \circ V - \bw_B \circ f(X, U)  = & 0  \\
& \bw_B \circ\Omega + \bw_B \circ\, \partial_x f(X, U) \circ \Lambda = & 0 \\
& \lambda^b - \lambda^a - \bw_B^T \Omega  = & 0\\
& x^b - x^a - \bw_B^T V  = & 0\\
&\bw_B \circ\,\partial_u f(X, U)\circ \Lambda  = & 0\\
& e(x^a, x^b) \le & 0\\
& \lambda^a +  \partial_{x^a} \overline{E}(\nu, x^a, x^b) = &0  \\
& \lambda^b - \partial_{x^b} \overline{E}(\nu, x^a, x^b) = & 0\\
&  \nu\ \dagger  e(x(\tau^a), x(\tau^b)) &
\end{array} \right. & \label{eq:probPlamN_a*b*}
\end{eqnarray}
%
%----------------
\begin{remark}\label{rem:notation4PlamNs}
Given a generalized root-finding problem in the form of Problem~$(P^\lambda)$, its discretization over the grid $\pi^N$ is denoted by Problem~$\left({P_{a^*,b^*}^{\lambda, N}}\right)$, where the subscripts $a^*, b^*$ imply that the state trajectory is represented by the $a$-form of its spectral discretization (denoted by $a^*$) while the costate trajectory is represented by the $b$-form of its spectral discretization (denoted by $b^*$).
\end{remark}
%-----------------
\begin{remark}\label{rem:16PlamNs}
In much the same way as was noted in Remark~\ref{rem:fourPNs} that Problem~$(P)$ can be ``directly'' discretized in at least four different ways via the spectral and pseudospectral methods associated with the $a$- and $b$-forms of the Birkhoff interpolants, it follows from Remark~\ref{rem:notation4PlamNs} that an ``indirect'' Birkhoff method can be generated in at least sixteen different ways.  In using the notation described in Remark~\ref{rem:notation4PlamNs}, these sixteen discretizations can be easily derived and represented as Problems~$\left({P_{a,a}^{\lambda, N}}\right), \left({P_{a,b}^{\lambda, N}}\right), \left({P_{b,a}^{\lambda, N}}\right), \left({P_{b,b}^{\lambda, N}}\right), \left({P_{a^*,a^*}^{\lambda, N}}\right), \left({P_{a^*,b^*}^{\lambda, N}}\right)$, \\
$\left({P_{b^*,a^*}^{\lambda, N}}\right), \left({P_{b^*,b^*}^{\lambda, N}}\right), \left({P_{a^*,a}^{\lambda, N}}\right)$,
$\left({P_{a^*,b}^{\lambda, N}}\right), \left({P_{b^*,a}^{\lambda, N}}\right), \left({P_{b^*,b}^{\lambda, N}}\right),
\left({P_{a,a^*}^{\lambda, N}}\right)$, \\
$\left({P_{a,b^*}^{\lambda, N}}\right), \left({P_{b,a^*}^{\lambda, N}}\right)$ and  $\left({P_{b,b^*}^{\lambda, N}}\right)$
\end{remark}
%-----------------
%

From Remarks~\ref{rem:fourPNs} and \ref{rem:16PlamNs} it follows that given Problem~$(P)$, one can construct at least twenty different Birkhoff discretizations over the same grid $\pi^N$.  In the next section, we show that these discretizations are related to each other in a specific way via the CMP. As outlined in Section~\ref{sec:intro}, the spectral algorithm to solve Problem~$(P)$ relies heavily on the CMP.

\section{The Covector Mapping Principle and Its Birkhoff Theorems}
\label{sec:CMP-derive}
If a discretization satisfies the CMP, then dualization and discretization commute\cite{ross-book, perspective, ross:roadmap-2005,CMP-CDC-2006}.  To explain, in simpler terms, why commutative operations in trajectory optimization are crucial\cite{ross:roadmap-2005} to generate the correct solution, consider the following simple univariate function,
\begin{equation}
g^N(x) := \frac{N x + x^2}{N^2}, \quad x \in \Real, \ N \in \mathbb{N}
\end{equation}
Suppose we minimize $g^N(x)$ with respect to $x$ (holding $N$ constant).  It is straightforward to show that this yields,
\begin{equation}\label{eq:mingN}
\min_x g^N(x) = -\frac{1}{4}
\end{equation}
Obviously, taking the limit of the left hand side of \eqref{eq:mingN} as $N \to \infty$ generates,
\begin{align}
%\min_x g^N(x) &= -\frac{1}{4}  \\
 \lim_{N \to \infty}\left(\min_x g^N(x) \right) &= -\frac{1}{4}
\end{align}
%\end{subequations}
%
That is, minimizing $g^N(x)$ before taking the limit $N \to \infty$ produces the number, $-1/4$; however, if we reverse the operations, we get a different answer:
\begin{subequations}
\begin{align}
 \lim_{N \to \infty} g^N(x) &= 0 \\
 \min_x \left( \lim_{N \to \infty} g^N(x) \right) &= 0
\end{align}
\end{subequations}
Obviously,
\begin{equation}
\lim_{N \to \infty}\min_x g^N(x) \ne \min_x \lim_{N \to \infty} g^N(x)
\end{equation}
That is, $\lim$ and $\min$ are not commutative operations.  Note that $\arg\min$ and $\lim$ are also not commutative:
\begin{equation}
\arg\min_x \left(\lim_{N \to \infty} g^N(x)\right) \ne \lim_{N \to \infty} \left(\arg\min_x g^N(x)\right)
\end{equation}
Trajectory optimization involves a number of commutative operations\cite{ross-book,ross:roadmap-2005,perspective}; hence, it is crucial that any proposed trajectory optimization method satisfy the CMP. Furthermore, the guess-free spectral algorithm\cite{spec-alg,DIDO:arXiv,fastmesh,furtherResults,ross:guess-free} for solving a generic optimal control problem rests on an explicit computational equivalence between direct and indirect methods.  In other words, the CMP is an integral part of a modern spectral algorithm\cite{DIDO:arXiv}. To develop such equivalence theorems, we follow the standard procedure\cite{ross-book,perspective, PSReview-ARC-2012,ross:roadmap-2005,CMP-CDC-2006,arb-grid} in using the CMP. An application of this procedure for the construction of a Birkhoff theory for optimal control can be cast in terms of the following steps:
%---------------
\begin{enumerate}
\item Develop Problem~$(P^\lambda)$ that results from an application of Pontryagin's principle to Problem~$(P)$.
\item Discretize Problem~$(P^\lambda)$ generated in Step~1 by some Birkhoff method.
\item Develop the stationarity conditions for Problems~$(P^N_\theta), \, \theta \in \set{a, b, a^*, b^*}$ by applying the Karush-Kuhn-Tucker (KKT) Theorem (or alternatively, the Fritz John (FJ) theorem).
\item Deduce the conditions that generate a relationship between the KKT (or FJ) multipliers of Step~3 and the Birkhoff-discretized covectors used in Step~2.
\end{enumerate}
%----------------
Steps~1 and 2 are presented in Section~\ref{sec:16PlamNs}.  In this section, we develop Steps~3 and 4. For the purpose of brevity we develop just three equivalence theorems.  A total of at least sixteen theorems (cf.~Remark~\ref{rem:16PlamNs}) are possible.  All of these theorems require the following lemmas and theorem:
%=======================================
\begin{lemma}\label{lemma:initTVC}
Let, $\psi_a \in \Real$ be defined by
\begin{equation}
\psi_a := \psi_b - \bw_B^T \Psi_B
\end{equation}
Then, $\psi_a =  \psi_b - \Psi_B^T W_B \bb$.
\end{lemma}
%-----------------
\begin{proof}
The proof of this lemma follows quite simply from the fact that $W_B$ is diagonal and $\bb$ is a vector of ones.
\end{proof}
%=======================
\begin{lemma}\label{lemma:int-by-parts=approx}
Let $\bdelta^N$ be defined as in Lemma~\ref{lemma:pre-int-by-parts}.  Let Hypothesis~\ref{hyp:piN=quadPts} hold. Then, given any $\epsilon > 0$, there exists an $N_\epsilon \in \mathbb{N} $ such that $\norm{\bdelta^{N}}_\infty \le \epsilon$ for all  $N \ge N_\epsilon $.
\end{lemma}
%------------------
\begin{proof}
The result follows directly from Hypothesis~\ref{hyp:piN=quadPts} and Lemma~\ref{lemma:pre-int-by-parts}.
\end{proof}
%==================
\begin{theorem}\label{theorem:ross-bigT}
Let Hypothesis 2 hold. Then, given any $\epsilon > 0$, there exists an $N_\epsilon \in \mathbb{N} $ such that for all  $N \ge N_\epsilon $,
\begin{equation}\label{eq:int-by-parts=exact}
W_B \bB^b + [\bB^a]^T W_B = 0
\end{equation}
\end{theorem}
%------------------
\begin{proof}
This theorem is a direct consequence of Lemma~\ref{lemma:int-by-parts=approx} and Definition~\ref{def:approx=exact}.
\end{proof}
%=====================================

\subsection{Spectral Covector Mapping Theorems for Problem~$(P^N_a)$}
In \cite{advances}, a number of key results were presented in using weighted inner products for constructing a Lagrangian relevant to a unified development of Lagrange PS methods (i.e., PS methods based on differentiation matrices). In using the insights of \cite{advances}, we formulate the following specially weighted Lagrangian for Problem~$(P^N_a)$:
\begin{multline}\label{eq:Lagrangian-WB}
L(\Psi_B, \Psi_V, \psi_b, \psi_e, X, U, V, x^a, x^b) := E(x^a, x^b) \\
+ \Psi_B^T W_B \left(X - \bB^aV - x^a\bb \right)
+ \Psi_V^T W_B \left(f(X,U) - V\right) \\
+ \psi_b\left(x^a + \bw_B^T V - x^b \right) + \psi_e e(x^a, x^b)
\end{multline}
where, $\Psi_B \in \real{N+1}, \Psi_V \in \real{N+1}, \psi_b \in \Real $ and  $\psi_e \in \Real$ are the Lagrange multipliers associated with the appropriate constraint equations implied in \eqref{eq:Lagrangian-WB}.  Furthermore, $\psi_e$ satisfies the complementarity condition, $\psi_e \dagger  e(x^a, x^b)$.  Taking the derivative of this Lagrangian with respect to the primal variables and setting it to zero generates the stationarity conditions.  Carrying out this exercise, it is straightforward to derive the following equations:
\begin{subequations}\label{eq:gradLwts=0}
\begin{align}
\partial_XL &= W_B\Psi_B + \partial_x f(X, U) \circ W_B \Psi_V  =  0\\
\partial_UL &= \partial_u f(X, U)\circ W_B \Psi_V = 0 \\
\partial_VL & = -\big[\bB^a\big]^T W_B \Psi_B - W_B \Psi_V + \psi_b \bw_B = 0 \label{eq:gradLVwt}\\
\partial_{x^a}L & = \partial_{x^a} \overline{E}(\psi_e, x^a, x^b) - \Psi_B^T W_B \bb + \psi_b = 0 \label{eq:initTVCviaKKT-1}\\
\partial_{x^b}L & = \partial_{x^b} \overline{E}(\psi_e, x^a, x^b) - \psi_b  = 0 \label{eq:finalTVCviaKKT-1}
\end{align}
\end{subequations}
%

%=================
\begin{theorem}[Birkhoff Pseudospectral-Spectral Covec. Map. Thm.] \label{thm:CMT-1}
Assume the following:
%-----------------
\begin{enumerate}
\item Hypotheses~\ref{hyp:a=b} and \ref{hyp:piN=quadPts}  hold.
\item Problem~$(P)$ is discretized according to Problem~$(P^N_a)$.
\item The Lagrangian given by \eqref{eq:Lagrangian-WB} holds over a finite-dimensional inner-product space defined by the Birkhoff quadrature weights.
\end{enumerate}
%-----------------
Then, there exists an $N_\epsilon \in \mathbb{N}$ such that for all $N \ge N_\epsilon$, the first-order optimality conditions for Problem~$(P^N_a)$ are identical to the pseudospectral-spectral approximation of Problem~$(P^\lambda)$ given by  Problem~$(P^{\lambda,N}_{a,b^*})$  under  the following equivalency conditions:
\begin{subequations}\label{eq:CMT-1}
\begin{align}
\psi_a &= \lambda^a   & \psi_b &= \lambda^b  & \psi_e &= \nu \\
\psi_V & = \Lambda    & \psi_B & = \Omega
\end{align}
\end{subequations}
\end{theorem}
%=====================
\begin{proof}
From the results of Section~\ref{sec:16PlamNs} and Remark~\ref{rem:16PlamNs}, it follows that Problem~$(P^{\lambda,N}_{a,b^*})$ can be written as,
%--------------------------------------
\begin{eqnarray}
& (\textsf{$P_{a, b^*}^{\lambda, N}$}) \left\{
\begin{array}{lrl}
&x^a\,\bb + \bB^a V - X = &0 \\
&\lambda^b \bw_B \circ \bb + \bw_B \circ \left(\bB^b \Omega \right)  - \bw_B \circ \Lambda= & 0  \\
& f(X, U)- V  = & 0 \\
& \bw_B \circ \Omega + \bw_B \circ\, \partial_x f(X, U) \circ \Lambda = & 0 \\
& \lambda^b - \lambda^a - \bw_B^T \Omega  = & 0\\
& x^b - x^a - \bw_B^T V  = & 0\\
&\bw_B \circ \,\partial_u f(X, U)\circ \Lambda  = & 0\\
&  e(x^a, x^b) \le & 0\\
& \lambda^a + \partial_{x^a} \overline{E}(\nu, x^a, x^b) = & 0\\
& \lambda^b - \partial_{x^b} \overline{E}(\nu, x^a, x^b) = & 0\\
& \nu\ \dagger  e(x(\tau^a), x(\tau^b))
\end{array} \right. & \label{eq:probPlamN_ab*}
\end{eqnarray}
%
%----------------
The first-order optimality conditions for Problem~$(P^N_a)$ are given by \eqref{eq:gradLwts=0} together with the complementarity condition, $\psi_e \dagger  e(x^a, x^b)$.  Hence, it suffices to show that the stationarity conditions in \eqref{eq:probPlamN_ab*} are the same as \eqref{eq:gradLwts=0} under the equivalency conditions given by \eqref{eq:CMT-1}.

Using Lemma~\ref{lemma:initTVC}, \eqref{eq:initTVCviaKKT-1} can be rewritten as,
\begin{equation}\label{eq:initTVCviaKKT-2}
\partial_{x^a} \overline{E}(\psi_e, x^a, x^b) + \psi_a = 0
\end{equation}
%
%Hence, from \eqref{eq:initTVCviaKKT-2} and \eqref{eq:finalTVCviaKKT-1} it follows that $\psi_a$ and $\psi_b$ satisfy the same conditions as $\lambda_a$ and $\lambda_b$ respectively.

Next, using \ref{eq:int-by-parts=exact}, \eqref{eq:gradLVwt} can be written as,
\begin{equation}\label{eq:gradLVwtXformed-1}
 W_B\, \bB^b\, \Psi_B - W_B \Psi_V + \psi_b \bw_B = 0
\end{equation}
Because $\bw_B$ is identically equal to $W_B \bb$, \eqref{eq:gradLVwtXformed-1} can be rewritten as,
\begin{multline}\label{eq:gradLVwtXformed-2}
W_B\, \bB^b\, \Psi_B - W_B \Psi_V + \psi_b  W_B \bb \\
=  \bw_B \circ \bB^b\, \Psi_B -  \bw_B \circ  \Psi_V + \psi_b   \bw_B \circ  \bb  = 0
\end{multline}
Replacing \eqref{eq:initTVCviaKKT-1} by \eqref{eq:initTVCviaKKT-2}, \eqref{eq:gradLVwt} by \eqref{eq:gradLVwtXformed-2} and using the fact that $W_B \Upsilon = \bw_B \circ \Upsilon$ for any $\Upsilon \in \real{N+1}$, it follows that \eqref{eq:gradLwts=0} generates the stationarity conditions in $(P^{\lambda,N}_{a,b^*})$ under the substitutions given by \eqref{eq:CMT-1}.
\end{proof}
%======================

It is apparent that \eqref{eq:Lagrangian-WB} is also the unweighted Lagrangian for Problem~$(P^N_{a^*})$; hence, we have the following corollary:
%------------------
\begin{corollary}[Birkhoff Spectral Covector Mapping Theorem]
Assume the following:
%-----------------
\begin{enumerate}
\item Hypotheses~\ref{hyp:a=b} and \ref{hyp:piN=quadPts}  hold.
\item Problem~$(P)$ is discretized according to Problem~$(P^N_{a^*})$.
\item An unweighted Lagrangian given by \eqref{eq:Lagrangian-WB} holds over a finite-dimensional inner-product space.
\end{enumerate}
%-----------------
Then, there exists an $N_\epsilon \in \mathbb{N}$ such that for all $N \ge N_\epsilon$, the first-order optimality conditions for Problem~$(P^N_{a^*})$ are identical to the spectral approximation of Problem~$(P^\lambda)$ given by Problem~$(P^{\lambda,N}_{a^*,b^*})$ under the equivalency conditions given by \eqref{eq:CMT-1}.
\end{corollary}
%---------------------

\subsection{A Pseudospectral Covector Mapping Theorem for Problem~$(P^N_a)$}
Suppose we scale the $X, V$ and $U$ variables in Problem~$(P^N_a)$ by $W_B$ such that the new optimization variables $\widetilde{X}, \widetilde{V}$ and $\widetilde{U}$ are given by,
\begin{multline}\label{eq:scale-ross}
\widetilde{X} = \bw_B \circ X = W_B X, \qquad \widetilde{V} = \bw_B \circ V = W_B V, \\
\widetilde{U} = \bw_B \circ U = W_B U
\end{multline}
Now consider the specially weighted Lagrangian for the scaled Problem~$(P^N_a)$:
\begin{multline}\label{eq:Lagrangian-noWts}
\widetilde{L}(\widetilde{\Psi}_B, \widetilde{\Psi}_V, \widetilde{\psi}_b, \widetilde{\psi}_e, \widetilde{X}, \widetilde{U}, \widetilde{V}, x^a, x^b)\\
 := E(x^a, x^b)
+ \widetilde{\Psi}_B^T W_B \left(W_B^{-1} \widetilde{X} - \bB^a W_B^{-1} \widetilde{V} - x^a\bb \right)\\
+ \widetilde{\Psi}_V^T  \left(f(W_B^{-1} \widetilde{X},\, W_B^{-1} \widetilde{U}) - W_B^{-1} \widetilde{V}\right) + \widetilde{\psi}_b\left(x^a + \bw_B^T W_B^{-1} \widetilde{V} - x^b \right) \\
+ \widetilde{\psi}_e e(x^a, x^b)
\end{multline}
where, $\widetilde{\Psi}_B \in \real{N+1}, \widetilde{\Psi}_V \in \real{N+1}, \widetilde{\psi}_b \in \Real $ and  $\widetilde{\psi}_e \in \Real$ are the Lagrange multipliers associated with the appropriate constraint equations implied in \eqref{eq:Lagrangian-noWts} with $\widetilde{\psi}_e$ satisfying the complementarity condition, $\widetilde{\psi}_e \dagger  e(x^a, x^b)$. Taking the derivative of this Lagrangian with respect to the scaled primal variables and setting it to zero generates the stationarity conditions that are analogous to \eqref{eq:gradLwts=0}.  These equations are given by,
\begin{subequations}\label{eq:gradL=0}
\begin{align}
\partial_{\widetilde{X}}\widetilde{L} &=  \widetilde{\Psi}_B + \partial_x f(X, U) \circ  W_B^{-1} \widetilde{\Psi}_V  =  0\\
\partial_{\widetilde{U}} \widetilde{L} &= \partial_u f(X, U)\circ  W_B^{-1} \widetilde{\Psi}_V = 0 \\
\partial_{\widetilde{V}} \widetilde{L} & = - W_B^{-1} \big[\bB^a\big]^T W_B \widetilde{\Psi}_B - W_B^{-1} \widetilde{\Psi}_V + \widetilde{\psi}_b\, \bb = 0 \label{eq:gradLV}\\
\partial_{x^a} \widetilde{L} & = \partial_{x^a} \overline{E}(\widetilde{\psi}_e, x^a, x^b) - \widetilde{\Psi}_B^T \, \bw_B + \widetilde{\psi}_b = 0 \label{eq:initTVCviaKKT-21}\\
\partial_{x^b} \widetilde{L} & = \partial_{x^b} \overline{E}(\widetilde{\psi}_e, x^a, x^b) - \widetilde{\psi}_b  = 0 \label{eq:finalTVCviaKKT-2}
\end{align}
\end{subequations}
where, $\widetilde{L}$ is an abbreviation for  $\widetilde{L}(\widetilde{\Psi}_B, \widetilde{\Psi}_V, \widetilde{\psi}_b, \widetilde{\psi}_e, \widetilde{X}, \widetilde{U}, \widetilde{V}, x^a, x^b)$

%=================
\begin{theorem}[Birkhoff Pseudospectral Covector Mapping Theorem]
Assume the following:
%-----------------
\begin{enumerate}
\item Hypotheses~\ref{hyp:a=b} and \ref{hyp:piN=quadPts}  hold.
\item Problem~$(P)$ is discretized according to Problem~$(P^N_a)$.
\item The optimization variables of Problem~$(P^N_a)$ are scaled according to \eqref{eq:scale-ross}.
\item The Lagrangian given by \eqref{eq:Lagrangian-noWts} holds over a finite-dimensional inner-product space.
\end{enumerate}
%-----------------
Then, there exists an $N_\epsilon \in \mathbb{N}$ such that for all $N \ge N_\epsilon$, the first-order optimality conditions for (the variable-scaled) Problem~$(P^N_a)$ are identical to the pseudospectral approximation of Problem~$(P^\lambda)$ given by Problem~$(P^{\lambda,N}_{a,b})$ under the following equivalency conditions:
\begin{subequations}\label{eq:CMT-2}
\begin{align}
\widetilde{\psi}_a &= \lambda^a   & \widetilde{\psi}_b &= \lambda^b  & \widetilde{\psi}_e &= \nu \\
\widetilde{\psi}_V & = \bw_B \circ \Lambda    & \widetilde{\psi}_B & =  \Omega
\end{align}
\end{subequations}
\end{theorem}
%=====================
\begin{proof}
Following Remark~\ref{rem:16PlamNs}, Problem~$(P^{\lambda,N}_{a,b})$ can be written as,
%--------------------------------------
%
\begin{eqnarray}
& (\textsf{$P_{a,b}^{\lambda, N}$}) \left\{
\begin{array}{lrl}
& X  = &x^a\,\bb + \bB^a V \\
& \Lambda = &\lambda^b\, \bb + \bB^b \Omega \\
& V =& f(X, U)  \\
& -\Omega = & \partial_x f(X, U)\circ \Lambda\\
& \lambda^b = & \lambda^a + \bw_B^T \Omega \\
& x^b = & x^a + \bw_B^T V \\
& 0 = &\partial_u f(X, U)\circ \Lambda \\
& 0 \ge &  e(x^a, x^b)\\
& -\lambda^a = & \partial_{x^a} \overline{E}(\nu, x^a, x^b) \\
& \lambda^b = & \partial_{x^b} \overline{E}(\nu, x^a, x^b) \\
& & \nu \dagger  e(x(\tau^a), x(\tau^b))
\end{array} \right. & \label{eq:probPlamN_ab}
\end{eqnarray}
The steps in the remainder of the proof are quite similar to that of Theorem~\ref{thm:CMT-1} and hence omitted for brevity.
\end{proof}
%=======================

%---------------
\begin{remark}
Equations~\eqref{eq:probPlamN_ab} and \eqref{eq:probPlamN_ab*} may be viewed as a more granular case of DIDO's generalized equation developed in \cite{DIDO:arXiv}.
\end{remark}
%-------------

\begin{figure}[h!]
      \centering
      {\parbox{0.9\columnwidth}{
      \centering
      {\includegraphics[width = 0.95\columnwidth]{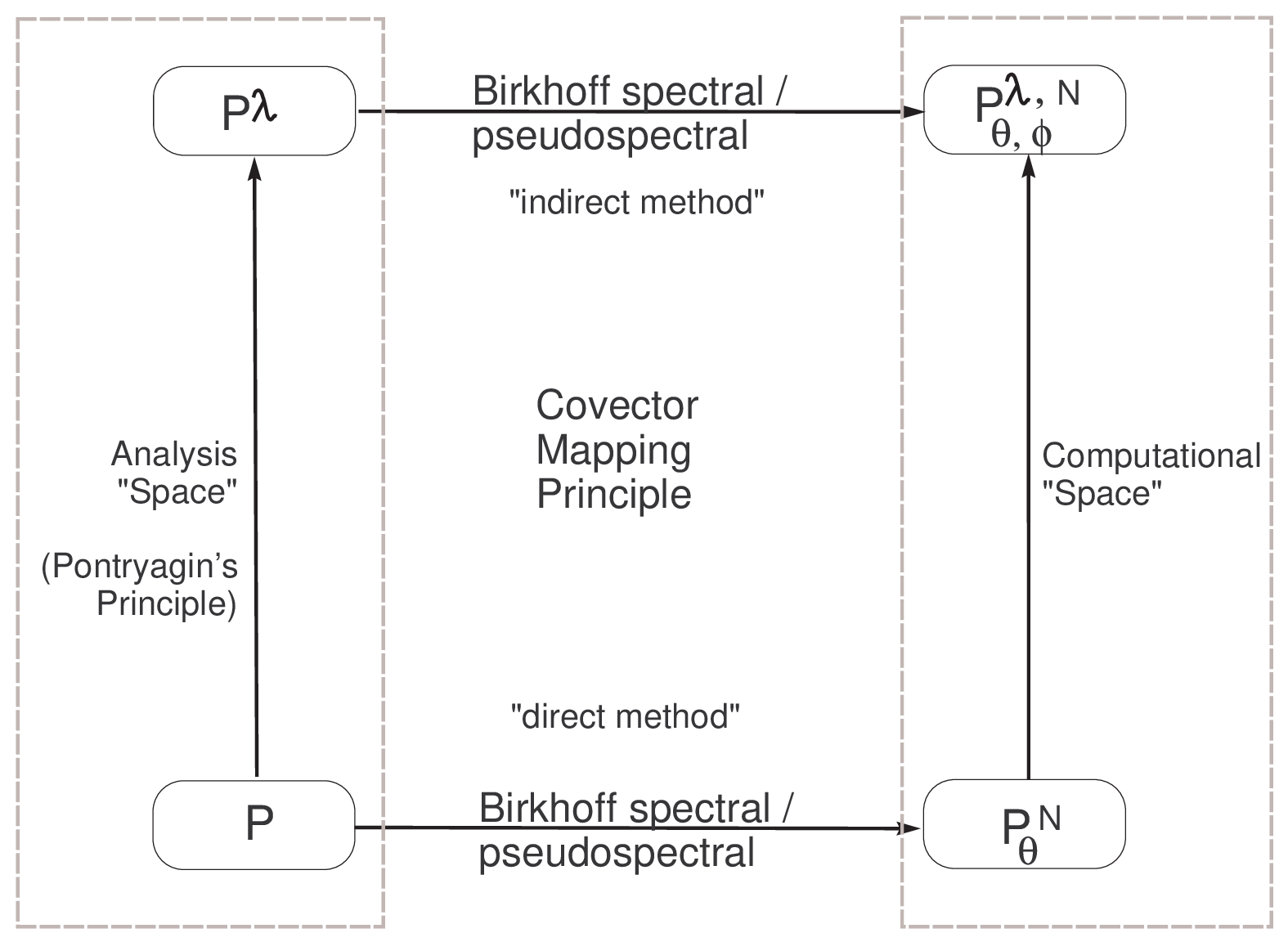}}
      \caption{{Schematic for the CMP illustrating how discretization commutes with dualization for Birkhoff methods with $\theta,\phi \in \set{a, b, a^*, b^*}$. Figure adapted from \cite{ross-book} and \cite{DIDO:arXiv} with permission of I. M. Ross \copyright\ I. M. Ross, 2015, 2020.}}\label{fig:CMP4Birkhoff}
      }
      }
\end{figure}

On the basis of Remarks~\ref{rem:fourPNs} and \ref{rem:16PlamNs} it is apparent that at least thirteen more covector mapping theorems can be developed using different spectral and pseudospectral discretizations of Problems~$(P)$ and $(P^\lambda)$.  All of these theorems hold for any grid $\pi^N$ that satisfies Hypotheses~\ref{hyp:a=b} and \ref{hyp:piN=quadPts}.  The collection of these theorems is depicted pictorially in Fig.~\ref{fig:CMP4Birkhoff} as the CMP.
%
%======================================================================================
%==========================================================================================
%

Once a valid choice of $\pi^N$ is made, then a selection of a particular covector mapping theorem can be extracted out of Fig.~\ref{fig:CMP4Birkhoff} to construct a guess-free fast spectral algorithm\cite{spec-alg,DIDO:arXiv,ross:guess-free,fastmesh,furtherResults} that is customized to a particular Birkhoff implementation.  Such implementation details are described in part II of this paper\cite{newBirk-part-II}.

\section{Conclusions}
Although pseudospectral (PS) theory based on Lagrange interpolants (i.e., PS theory based on differentiation matrices) are already used in mission design and flight operations, its practical implementation is steeped in various layers of sophistication that manage the growth of condition numbers and round-off errors. It has long been known that the source of these challenges is indeed the starting point; i.e., the use of Lagrange interpolants.  Lagrange interpolants also form the basis of Runge-Kutta methods, albeit in their low-order forms. Consequently, an escape from these challenges seemed improbable.
In choosing an altogether new starting point, the Birkhoff theory attacks the problem right at the source.
This paper shows that the recently proposed universal Birkhoff interpolant may indeed be used as a new starting point for trajectory optimization.
The new first-principles' approach motivates and integrates many number of disparate concepts in applied mathematics: Birkhoff interpolants, projections in function spaces, spectral (Galerkin) and pseudospectral approximation theories, Pontryagin's principle, Karush-Kuhn-Tucker theorem, weighted inner-products in finite-dimensional spaces and the covector mapping principle.  Despite using a broad swath of results across the mathematical spectrum, a covector mapping theorem can be stated quite simply in terms of its founding principle: a proper discretization can be commuted with dualization. Consequently, many elements of the previously-developed, guess-free, fast spectral algorithm can be easily adapted to the new framework.
More importantly, the equivalences between direct/indirect as well as spectral/pseudospectral methods provide additional avenues for algorithmic acceleration.  Combined with the isolation of the Birkhoff-specific computations to a well-conditioned linear system, a Birkhoff-theoretic method is poised to generate fast and accurate solutions to nonlinear, nonconvex trajectory optimization problems.

%%%%%%%%
%Bibliography
%%%%%%%%

\begin{spacing}{1.0}

\end{spacing}


\begin{thebibliography}{10}
\newcommand{\enquote}[1]{``#1''}

\bibitem{ross-book}
I. M. Ross, \textit{A Primer on Pontryagin's Principle in Optimal Control}, Second Edition, Collegiate Publishers, San Francisco, CA, 2015. https://www.amazon.com/Primer-Pontryagins-Principle-Optimal-Control/dp/0984357114.

\bibitem{longuski}
J. M. Longuski, J. J. Guzm\'{a}n and J. E. Prussing, \textit{Optimal Control with Aerospace Applications}, Springer, New York, N.Y., 2014.


\bibitem{trelat-2017}
R. Bonalli, B. H\'{e}riss\'{e} and E. Tr\'{e}lat, ``Analytical Initialization of a Continuation-Based Indirect Method for Optimal Control of Endo-Atmospheric Launch Vehicle Systems,'' \textit{IFAC-PapersOnLine}, Vol.50, No.~1, 2017, pp.~482-487. https://doi.org/10.48550/arXiv.1703.05117.


\bibitem{asteroid-trajOpt}
A. D'Ambrosio, A. Carbone, F. Curti, ``Optimal Maneuvers Around Binary Asteroids Using Particle Swarm Optimization and Machine Learning,'' \textit{Journal of Spacecraft and Rockets}, 2023, 60/5, pp.~1458--1472. https://doi.org/10.2514/1.A35317.


\bibitem{ddp-2021}
G. I. Boutselis and E. Theodorou, ``Discrete-Time Differential Dynamic Programming on Lie Groups: Derivation, Convergence Analysis, and Numerical Results,'' \textit{IEEE Transactions On Automatic Control}, Vol. 66, NO. 10, October 2021, pp.~4636-4651. http://dx.doi.org/10.1109/TAC.2020.3034206.


\bibitem{trelat:guidance}
R. Bonalli, B. H\'{e}riss\'{e} and E. Tr\'{e}lat, ``Optimal Control of Endoatmospheric Launch Vehicle Systems: Geometric and Computational Issues,'' \textit{IEEE Transactions on Automatic Control}, vol. 65, no. 6, pp.~2418-2433, June 2020. https://dx.doi.org/10.1109/tac.2019.2929099.

\bibitem{convex-2015}
S. J. Zhang, B. A\c{c}ikmes\c{s}e, S. S.-M. Swei and D. Prabhu, ``Convex programming approach to real-time trajectory optimization for mars aerocapture,'' \textit{2015 IEEE Aerospace Conference}, Big Sky, MT, USA, 2015, pp. 1-7. http://dx.doi.org/10.1109/AERO.2015.7119111


\bibitem{conway:survey}
B. A. Conway, ``A Survey of Methods Available for the Numerical Optimization of Continuous Dynamic Systems,'' \textit{Journal of Optimization Theory and Applications}, Vol.~152, 2012, pp.~271--306. https://doi.org/10.1007/s10957-011-9918-z.


\bibitem{trelat:survey}
E. Tr\'{e}lat, ``Optimal Control and Applications to Aerospace: Some Results and Challenges,'' \textit{Journal
of Optimization Theory and Applications}, 2012, 154 (3), pp~713-758. DOI:10.1007/s10957-012-0050-5


\bibitem{flatness-2013}
D. Desiderio and M. Lovera, ``Guidance and Control for Planetary Landing: Flatness-Based Approach,'' \textit{IEEE Transactions on Control Systems Technology}, vol. 21, no. 4, pp. 1280--1294, July 2013. http://dx.doi.org/10.1109/TCST.2012.2202664


\bibitem{convex-2011}
B. A\c{c}ikmes\c{s}e and L. Blackmore, ``Lossless convexification of a class of optimal control problems with non-convex control constraints,'' \textit{Automatica}, vol. 47, no. 2, 2011, pp. 341--347. http://dx.doi.org/10.1109/ACC.2012.6314722



\bibitem{perspective}
I. M. Ross and F. Fahroo, ``A Perspective on Methods for Trajectory Optimization,'' \textit{AIAA/AAS Astrodynamics Specialist Conference and Exhibit}, 5-8 August, 2002, Monterey, CA. AIAA 2002-4727. https://doi.org/10.2514/6.2002-4727.



\bibitem{seywald-94}
H. Seywald, ``Trajectory  optimization  based  on  differential  inclusion,'' \textit{Journal of Guidance, Control, and Dynamics}, Vol. 17, No. 3, 1994, pp.~480--487. https://doi.org/10.2514/3.21224

\bibitem{vonStryk-92}
O. Von Stryk, R. Bulirsch, ``Direct and indirect methods for trajectory optimization,'' \textit{Annals of operations research}, Vol. 37, No. 1, 1992, pp. 357--373. http://dx.doi.org/10.1007/BF02071065



\bibitem{lu:entry-guidance}
Lu, P., ``Entry Guidance: A Unified Method,'' \textit{Journal of Guidance, Control, and Dynamics}, Vol. 37, No.3, 2014, pp.~713--728. http://dx.doi.org/10.2514/1.62605


\bibitem{lu:editorial}
Lu, P., ``What is Guidance,'' \textit{Journal of Guidance, Control, and Dynamics}, Vol. 44, No.7, July 2021, pp.~1237--1238. https://doi.org/10.2514/1.G006191

\bibitem{RTOC-jgcd}
I. M. Ross, P. Sekhavat, A. Fleming and Q. Gong, Q., ``Optimal Feedback Control: Foundations, Examples and Experimental Results for a New Approach,'' \textit{Journal of Guidance, Control, and Dynamics}, Vol. 31, No.2, Mar-Apr 2008, pp. 307-321. http://dx.doi.org/10.2514/1.29532


\bibitem{HJB=PP=DIDO}
J. Weston, D. Toli\'{c}, D. and I. Palunko, ``Application of Hamilton-Jacobi-Bellman Equation/Pontryagin's Principle for Constrained Optimal Control,'' \textit{Journal of  Optimization Theory and Applications}, 200, 2024, pp.~437--462. http://dx.doi.org/10.1007/s10957-023-02364-4


\bibitem{mease:guidance-compute}
F. Najson, K. D. Mease, ``Computationally Inexpensive Guidance Algorithm for Fuel-Efficient Terminal Descent,'' \textit{Journal of Guidance, Control, and Dynamics}, Vol. 29, No.4, 2006, pp.~955-964.
http://dx.doi.org/10.2514/6.2005-6289


\bibitem{brysonHo}
Bryson, A. E., and Ho, Y.-C., \textit{Applied Optimal Control}, Hemisphere, New York, 1975 (Revised
Printing; original publication, 1969).

\bibitem{clarke-2013book}
F. Clarke, \textit{Functional Analysis, Calculus of Variations and Optimal Control}, Springer-Verlag, London, 2013.

\bibitem{vinter}
R. B. Vinter,  \textit{Optimal Control}, Birkh\"{a}user, Boston, MA, 2000.

\bibitem{LGR-rtoc}
I. M. Ross, Q. Gong, F. Fahroo and W. Kang, ``Practical stabilization through real-time optimal control,''  \textit{Proceedings of the 2006 American Control Conference}, Inst. of Electrical and Electronics Engineers, Piscataway, NJ, June 2006, pp.~14--16. http://dx.doi.org/10.1109/ACC.2006.1655372



\bibitem{PEG:history}
J. L. Goodman, ``Roland Jaggers and the Development of Space Shuttle Powered Explicit Guidance (PEG),'' AIAA SciTech Forum Virtual Event, January 11-15 \& 19-21, 2021. http://dx.doi.org/10.2514/6.2021-2021


\bibitem{lu:prop-guidance}
Lu, P., ``Propellant-Optimal Powered Descent Guidance,'' \textit{Journal of Guidance, Control, and Dynamics}, Vol. 41, No.4, Mar-Apr 2018, pp.~813--826. http://dx.doi.org/10.2514/1.G007214



\bibitem{bollinoRossDoman-2006}
K. P. Bollino, I. M. Ross and D. Doman, ``Optimal nonlinear feedback guidance for reentry vehicles,''
\textit{AIAA Guidance, Navigation, and Control Conference and Exhibit}, Keystone, CO, AIAA-2006-6074, 2006. http://dx.doi.org/10.2514/6.2006-6074

%\bibitem{bollino:PhD-2006}
%K. P. Bollino, \textit{High-Fidelity Real-Time Trajectory Optimization for Reusable Launch Vehicles}, Ph.D. Dissertation, Naval Postgraduate School, Monterey, CA, Dec 2006.

\bibitem{bollinoRoss-2007}
K. P. Bollino and I. M. Ross, ``A pseudospectral feedback method for real-time optimal guidance of reentry vehicles,'' \textit{Proc. of the 2007 IEEE American Control Conference}, New York, NY, Jul 2007. https://doi.org/10.1109/ACC.2007.4282500

\bibitem{lewis-NATO-2007}
L.R. Lewis and I.M. Ross, ``A pseudospectral method for real-time motion planning and obstacle avoidance,'' AVT-SCI Joint Symposium on Platform Innovations and System Integration for Unmanned Air, Land and Sea Vehicles, 2007. https://apps.dtic.mil/sti/pdfs/ADA478686.pdf


\bibitem{hurniSekRoss-2009}
Hurni, M.A., Sekhavat, P., Ross, I.M. (2010). An Info-Centric Trajectory Planner for Unmanned Ground Vehicles. In: Hirsch, M., Pardalos, P., Murphey, R. (eds) Dynamics of Information Systems. Springer Optimization and Its Applications, vol 40. Springer, New York, NY. https://doi.org/10.1007/978-1-4419-5689-\_11.


%.A. Hurni, P. Sekhavat, and I.M. Ross, ``Pseudospectral optimal control algorithm for real-time trajectory planning,'' Proceedings of the 19th AAS/AIAA Space Flight Mechanics Meeting, 2009.

\bibitem{bollinoLewis-2008}
K.P. Bollino and L.R. Lewis, ``Collision-free Multi-UAV Optimal Path Planning and
Cooperative Control for Tactical Applications,'' \textit{Proceedings of the AIAA Guidance, Navigation and Control Conference and Exhibit}, AIAA 2008-7134, 2008. http://dx.doi.org/10.2514/6.2008-7134


\bibitem{bollinoLewis-2007}
K.P. Bollino and L.R. Lewis, ``Optimal path planning and control of tactical unmanned aerial vehicles in urban environments,'' \textit{Proceedings of the AUVSI's Unmanned Systems North America Conference}, 2007


\bibitem{bollinoEtAl-2007}
K. P. Bollino, L. R. Lewis, P. Sekhavat and I. M. Ross, ``Pseudospectral optimal control: a clear
road for autonomous intelligent path planning,'' \textit{AIAA Infotech@Aerospace 2007 Conference
and Exhibit}, Rohnert Park, CA, AIAA-2007-2831 (2007) http://dx.doi.org/10.2514/6.2007-2831

\bibitem{PSReview-ARC-2012}
I. M. Ross and M. Karpenko, ``A Review of Pseudospectral Optimal Control: From Theory to Flight,'' \textit{Annual Reviews in Control}, Vol.36, No.2, pp.182--197, 2012. http://dx.doi.org/10.1016/j.arcontrol.2012.09.002


\bibitem{spec-alg}
Q. Gong, F. Fahroo and I. M. Ross, ``Spectral Algorithm for Pseudospectral Methods in Optimal Control,'' \textit{Journal of Guidance, Control, and Dynamics}, vol. 31 no. 3, pp. 460-471, 2008. http://dx.doi.org/10.2514/1.32908

\bibitem{DIDO:arXiv}
I. M. Ross, ``Enhancements to the DIDO Optimal Control Toolbox,'' arXiv preprint, arXiv:2004.13112, 2020, https://arxiv.org/abs/2004.13112.


\bibitem{ross:roadmap-2005}
I. M. Ross, ``A Roadmap for Optimal Control: The Right Way to Commute,'' \textit{Annals of the New York Academy of Sciences}, 1065/1, 2005, 210--231. https://doi.org/10.1196/annals.1370.015


\bibitem{ross:guess-free}
I. M. Ross and Q. Gong, ``Guess-Free Trajectory Optimization,'' \textit{AIAA/AAS Astrodynamics Specialist Conference and Exhibit}, 18--21 August 2008, Honolulu, Hawaii, AIAA 2008-6273. https://core.ac.uk/download/pdf/36722475.pdf

\bibitem{Kang_2008_convergence}
W. Kang, I. M. Ross and Q. Gong, ``Pseudospectral Optimal Control and its Convergence Theorems,'' \textit{Analysis and Design of Nonlinear Control Systems}, Springer-Verlag, Berlin Heidelberg, 2008, pp.~109--126. http://dx.doi.org/10.1007/978-3-540-74358-3\_8


\bibitem{kang-rate-2010}
W. Kang, ``Rate of Convergence for a Legendre Pseudospectral Optimal Control of Feedback Linearizable Systems,'' \textit{Journal of Control Theory and Applications}, Vol.~8, No.~4, pp.~391--405, 2010. http://dx.doi.org/10.1007/s11768-010-9104-0


\bibitem{bellman-low-t}
I. M. Ross, Q. Gong and P. Sekhavat, ``Low-Thrust, High-Accuracy Trajectory Optimization,'' \textit{Journal of Guidance, Control and Dynamics}, Vol.~30, No.~4, pp.~921--933, 2007. http://dx.doi.org/10.2514/1.23181

\bibitem{bellman-conf}
I. M. Ross, Q. Gong and P. Sekhavat, ``The Bellman Pseudospectral Method,'' \textit{AIAA/AAS Astrodynamics Specialist Conference and Exhibit}, Honolulu, Hawaii, AIAA-2008-6448, August 18-21, 2008. https://doi.org/10.2514/6.2008-6448


\bibitem{TEI-JGCD-2011}
H. Yan, Q. Gong, C. Park, I. M. Ross, and C. N. D'Souza, ``High Accuracy Trajectory Optimization for a Trans-Earth Lunar Mission,'' \textit{Journal of Guidance, Control and Dynamics}, Vol. 34, No. 4, 2011, pp. 1219-1227. http://dx.doi.org/10.2514/1.49237


\bibitem{twist}
R. Baltensperger and M. R. Trummer, ``Spectral Differencing with a Twist,'' SIAM J. Sci. Comput., 24(5), 14651487, 2003. http://dx.doi.org/10.1137/S1064827501388182

\bibitem{trefethen-2000}
L. N. Trefethen, \textit{Spectral Methods in MATLAB}, SIAM, Philadelphia, PA, 2000.



\bibitem{rossJCAM-1}
I. M. Ross, ``An Optimal Control Theory for Nonlinear Optimization,'' \textit{Journal of Computational and Applied Mathematics}, Vol.~354, July 2019, 39--51. http://dx.doi.org/10.1016/j.cam.2018.12.044

\bibitem{trefethen-bau-1997}
L. N. Trefethen and D. Bau, III, \textit{Numerical Linear Algebra}, SIAM, Philadelphia, PA, 1997.


\bibitem{hesthavan}
J. S. Hesthaven, ``Integration Preconditioning Of Pseudospectral Operators. I. Basic Linear Operators,'' \textit{SIAM Journal of Numerical Analysis}, Vol.~35, No.~4, pp.~1571-1593, 1998. http://dx.doi.org/10.1137/S0036142997319182


\bibitem{elbarbary}
E. M. E. Elbarbary, ``Integration Preconditioning Matrix for Ultraspherical Pseudospectral Operators,'' \textit{SIAM Journal of Scientific Computaton}, Vol.~28, No.~3, pp.~1186-1201, 2006. https://doi.org/10.1137/050630982


\bibitem{wang}
L.-L Wang, M. D. Samson and X. Zhao, ``A Well-Conditioned Collocation Method Using a Pseudospectral Integration Matrix,'' \textit{SIAM Journal of Scientific Computaton}, Vol.~36, No.~3, pp.~A907-A929, 2014. http://dx.doi.org/10.1137/130922409

\bibitem{fornberg}
B. Fornberg,  A Practical Guide to Pseudospectral Methods. Cambridge University Press Cambridge, 1996.

\bibitem{boyd}
J. Boyd,  {\it Chebyshev and Fourier Spectral Methods,} Dover
Publications, Inc., Minola, New York, 2001.



\bibitem{lorentz}
G. G. Lorentz and K. L. Zeller, ``Birkhoff Interpolation,'' \textit{SIAM Journal of Numerical Analysis}, Vol.~8, No.~1, pp.~43-48, 1971. https://doi.org/10.1137/0708006

\bibitem{schoenberg}
I. J. Schoenberg, ``On Hermite-Birkhoff Interpolation,'' \textit{Journal of Mathematical Analysis and Applications}, Vol.~16, No.~3, pp.~538-543, 1966. https://doi.org/10.1016/0022-247X(66)90160-0


\bibitem{fastmesh}
N. Koeppen, I. M. Ross, L. C. Wilcox and R. J. Proulx, ``Fast Mesh Refinement in Pseudospectral Optimal Control,'' \textit{Journal of Guidance, Control, and Dynamics}, vol. 42 no. 4, pp. 711-722, 2018. http://dx.doi.org/10.2514/1.G003904

\bibitem{newBirk-2023}
I. M. Ross,  R. J. Proulx, C. F. Borges, ``A Universal Birkhoff Pseudospectral Method for Solving
Boundary Value Problems,'' \textit{Applied Mathematics and Computation}, Vol.~454, 128101, 2023. https://doi.org/10.1016/j.amc.2023.128101

\bibitem{furtherResults}
I. M. Ross and R. J. Proulx, ``Further Results on Fast Birkhoff Pseudospectral Optimal Control Programming,'' \textit{Journal of Guidance, Control, and Dynamics}, vol. 42 no. 9, pp. 2086--2092, 2019. http://dx.doi.org/10.2514/1.G004297



\bibitem{newBirk-part-II}
Proulx R. J. and Ross, I. M., ``Implementations of the Universal Birkhoff Theory for Fast Trajectory Optimization,'' \textit{Journal of Guidance, Control and Dynamics},  vol. 47, no. 12, pp. 2482--2496, Dec. 2024. https://doi.org/10.2514/1.G007738.



\bibitem{CMP-CDC-2006}
Q. Gong, I. M. Ross, W. Kang and F. Fahroo, ``On the Pseudospectral Covector Mapping Theorem for Nonlinear Optimal Control,'' \textit{Proceedings of the 45th IEEE Conference on Decision and Control}, 2006, pp. 2679-2686, doi: 10.1109/CDC.2006.377729.



\bibitem{scaling}
I. M. Ross, Q. Gong, M. Karpenko and R. J. Proulx, ``Scaling and Balancing for High-Performance Computation of Optimal Controls,''  \textit{Journal of Guidance, Control and Dynamics}, Vol.~41, No.~10, 2018, pp.~2086--2097. http://dx.doi.org/10.2514/1.G003382


%\bibitem{stryk-1992}
%von Stryk, O., Bulirsch, R., ``Direct and indirect methods for trajectory optimization,'' Annals of Opererations Research, 37, 1992, pp.~357--373.


\bibitem{cots-2017}
Cots, O., Gergaud, J.,  Goubinat, D., ``Direct and indirect methods in optimal control
with state constraints and the climbing trajectory of an aircraft,'' \textit{Optimal Control Applications and
Methods}, Wiley, 2018, 39, pp.~281--301. https://doi.org/10.1002/oca.2347

\bibitem{ross-CD}
Ross, I. M., ``Derivation of Coordinate Descent Algorithms from Optimal Control Theory'' \textit{Operartions Research Forum}, Vol.~4, No.~2, 2023. https://doi.org/10.1007/s43069-023-00215-6

\bibitem{TAC:linearizable}
Q. Gong, W. Kang and I. M. Ross, ``A Pseudospectral Method for
the Optimal Control of Constrained Feedback Linearizable Systems,''
\textit{IEEE Transactions on Automatic Control}, Vol.~51, No.~7,
July 2006, pp.~1115-1129. http://dx.doi.org/10.1109/TAC.2006.878570

\bibitem{paris-comparision-2006}
Paris, S. W., Riehl, J. P. and Sjaw, W. K., ``Enhanced Procedures for Direct Trajectory Optimization Using Nonlinear Programming and Implicit Integration,'' \textit{AIAA/AAS Astrodynamics Specialist Conference and Exhibit}, 21-24 August 2006, Keystone, CO, AIAA 2006-6309. http://dx.doi.org/10.2514/6.2006-6309

\bibitem{sandia-aas-23}
A. Javeed, D. Ridzal, D. Kouri and I. M. Ross, ``A Fast Matrix-Free Method for Low-Thrust Trajectory Optimization,'' \textit{AAS/AIAA Astrodynamics Specialist Conference and Exhibit}, Big Sky, MT, August 13-17, 2023, Paper No.~AAS-280. https://space-flight.org/docs/2023\_summer/Conference-Program.pdf

\bibitem{Bogaert-2014}
I. Bogaert, ``Iteration-Free Computation of Gauss--Legendre Quadrature Nodes and Weights,'' \textit{SIAM J. Sci. Comput.}, 36/3 (2014), A1008-A1026. https://doi.org/10.1137/140954969

\bibitem{olver-2021}
Olver, S., Slevinsky, R., and Townsend, A., ``Fast algorithms using orthogonal polynomials,'' \textit{Acta Numerica}, 29, 2021, pp.~573-699.
doi:10.1017/S0962492920000045

\bibitem{Gil-2019}
Gil, A., Segura, J. and Temme, N.M., ``Fast, reliable and unrestricted iterative computation of Gauss-Hermite and Gauss-Laguerre quadratures,'' \textit{Numerische Mathematik}, 143, 2019, pp.~649--682. https://doi.org/10.1007/s00211-019-01066-2

\bibitem{Hale-fast-comp-2013}
N. Hale, A. Townsend, Fast and Accurate Computation of Gauss--Legendre and Gauss--Jacobi Quadrature Nodes and Weights, SIAM J. Sci. Comput., 35/2 (2013), A652-A674. https://doi.org/10.1137/120889873

\bibitem{kreyszig-2011}
Kreyszig, E. \textit{Advanced Engineering Mathematics}, Wiley, 2011.

\bibitem{fft-original}
Cooley, J. W., and Tukey, J. W., ``An Algorithm for the Machine Calculation of Complex Fourier Series,'' \textit{Mathematics of Computation}, 19-90, 1965, 297--301. https://doi.org/10.2307/2003354

\bibitem{fftw}
Frigo, M. and Johnson, S. G., ``The Design and Implementation of FFTW3,'' \textit{Proceedings of the IEEE}, Vol.~93, No.~2, 2005,  216--231. https://doi.org/10.1109/JPROC.2004.840301

\bibitem{luenberger-2008}
Luenberger, D. G. and Ye, Y. \textit{Linear and Nonlinear Programming}, Springer, 2008.

\bibitem{NW:NumOptBook}
Nocedal, J. and Wright, S. \textit{Numerical Optimization},  Springer,  2006.

\bibitem{bazaraa-2006}
Bazaraa, M. S., Sherali, H. D. and Shetty, C. M. \textit{Nonlinear Programming: Theory and Algorithms}, Wiley-Interscience, 2006.

\bibitem{nesterov-book-2004}
Nesterov, Y., \textit{Introductory Lectures on Convex Optimization: A Basic Course}, Kluwer Academic Publishers, 2004.


\bibitem{advances}
F. Fahroo and I. M. Ross, ``Advances in Pseudospectral Methods for Optimal Control,'' \textit{AIAA Guidance, Navigation, and Control Conference}, AIAA Paper 2008-7309, Honolulu, Hawaii, August 2008. https://doi.org/10.2514/6.2008-7309

\bibitem{auto-knots}
Q. Gong and I. M. Ross, ``Autonomous Pseudospectral Knotting Methods for Space Mission Optimization,'' \textit{Advances in the Astronatuical Sciences}, Vol.~124, 2006, AAS 06-151, pp.~779--794.

\bibitem{Radau-GNC05}
F. Fahroo and I. M. Ross, ``Pseudospectral Methods for
Infinite-Horizon Nonlinear Optimal Control Problems,'' {\it Proceedings of the
AIAA Guidance, Navigation and Control Conference}, San Francisco, CA, August 15-18,
2005. https://arc.aiaa.org/doi/pdfplus/10.2514/6.2005-6076.

\bibitem{hnw-ode}
Hairer, E., N{\o}rsett, S. P. \& Wanner, G. \textit{Solving Ordinary Differential Equations I: Nonstiff Problems}. (Springer-Verlag, 1993).

\bibitem{ap-ode}
U. M. Ascher and L. R. Petzold, \textit{Computer Methods for Ordinary Differential Equations and Differential-Algebraic Equations}, SIAM, 1998.



\bibitem{cullum:1972}
J. Cullum, ``Finite-dimensional approximations of state constrainted
continuous optimal problems,'' SIAM J. Control, vol. 10, pp. 649--670,
1972. https://doi.org/10.1137/0310048

\bibitem{polak:book-1997}
E. Polak, \textit{Optimization: Algorithms and Consistent Approximations}, Heidelberg, Germany: Springer-Verlag, 1997.

\bibitem{boris:AMP:2004}
B. S. Mordukhovich and I. Shvartsman, ``The Approximate Maximum Principle in Constrained Optimal Control,'' SIAM Journal of Control and Optimization, Vol. 43, No. 3, 2004, pp. 1037--1062. https://doi.org/10.1137/S0363012903433012

\bibitem{CMP:history}
I. M. Ross, ``A Historical Introduction to the Covector Mapping Principle,'' \textit{Advances in the Astronautical Sciences}, Vol.~123, Univelt, San Diego, CA, 2006, pp.~1257--1278.

\bibitem{Gill:book}
P. E. Gill, W. Murray and M. H. Wright, \textit{Practical Optimization},
Academic Press, London, 1981.


\bibitem{arb-grid-aas}
Q. Gong, I. M. Ross and F. Fahroo, ``Pseudospectral Optimal Control On Arbitrary Grids,'' \textit{AAS Astrodynamics Specialist Conference}, AAS-09-405, 2009. https://core.ac.uk/download/pdf/36732678.pdf.

\bibitem{arb-grid}
Q. Gong, I. M. Ross and F. Fahroo, ``Spectral and Pseudospectral Optimal Control Over Arbitrary Grids,'' \textit{Journal of Optimization Theory and Applications}, vol.~169, no.~3, pp.~759-783, 2016. https://link.springer.com/article/10.1007/s10957-016-0909-y

\bibitem{nonpoly-interp}
I. H. Sloan, ``Nonpolynomial Interpolation,'' \textit{J. Approx. Theory}, 39, 1983, pp.~97--117. https://doi.org/10.1016/0021-9045(83)90085-0

\bibitem{rbf-interp}
T. A. Driscoll and B. Fornberg, ``Interpolation in the limit of increasingly flat radial basis functions,''
\textit{Computers \& Mathematics with Applications}, Volume 43, Issues 3-5, 2002, Pages 413--422. https://doi.org/10.1016/S0898-1221(01)00295-4

\bibitem{sinc-interp}
M. Sugihara and T. Matsuo, ``Recent developments of the Sinc numerical methods,''
\textit{Journal of Computational and Applied Mathematics}, Volumes 164-165, 2004, Pages 673--689. http://dx.doi.org/10.1016/j.cam.2003.09.016

\bibitem{wavelets}
Daubechies, I., \textit{Ten Lectures on Wavelets}, SIAM, 1992. doi:10.1137/1.9781611970104.

\bibitem{Radau-JGCD}
F. Fahroo and I. M. Ross, ``Pseudospectral Methods for Infinite-Horizon Optimal Control Problems,'' \textit{Journal of Guidance, Control and Dynamics}, Vol.~31, No.~4, pp.~927--936, 2008. http://dx.doi.org/10.2514/1.33117



\bibitem{knots}
I. M. Ross and F. Fahroo, ``Pseudospectral Knotting Methods for
Solving Optimal Control Problems,'' {\it Journal of Guidance,
Control and Dynamics}, Vol. 27, No. 3, pp. 397-405, 2004. https://doi.org/10.2514/1.3426

\bibitem{acc:hybrid}
I. M. Ross and F. Fahroo, ``Discrete Verification of Necessary
Conditions for Switched Nonlinear Optimal Control Systems,''
\textit{Proceedings of the American Control Conference}, June 2004,
Boston, MA. http://dx.doi.org/10.1109/ACC.2004.183014

\bibitem{hybrid:jgcd}
Ross, I. M., D'Souza, C. N., ``Hybrid Optimal Control Framework for Mission Planning,'' {\it Journal of Guidance, Control and Dynamics}, Vol. 28, No. 4, pp. 686--697, 2005. http://dx.doi.org/10.2514/1.8285



\bibitem{halkin:1966}
H. Halkin, ``A Maximum Principle of the Pontryagin Type for Systems
Described by Nonlinear Difference Equations,'' \textit{SIAM Journal of Control},
Vol. 4, No. 1, 1966, pp. 90--111. https://doi.org/10.1137/0304009






















\end{thebibliography}
\end{document}